\documentclass[a4paper]{article}
\usepackage{amsmath,amssymb,a4wide}
\usepackage{color,graphicx,epstopdf}
\usepackage[a4paper]{geometry}

\newcommand{\bu}{\boldsymbol u}
\newcommand{\bv}{\boldsymbol v}
\newcommand{\bw}{\boldsymbol w}

\newcommand{\be}{\boldsymbol e}

\newcommand{\bff}{\boldsymbol f}

\newcommand{\bphi}{\boldsymbol\phi}

\newtheorem{Theorem}{Theorem}
\newtheorem{lema}{Lemma}
\newcounter{remark}
\def\theremark {\arabic{remark}}
\newenvironment{remark}{\refstepcounter{remark}\par\noindent{\bf Remark\ \theremark}\ }{\par}
\newtheorem{Proof}{Proof}

\usepackage{hyperref}

\title{Error Analysis of Non Inf-sup Stable Discretizations of the
time-dependent Navier--Stokes Equations with Local Projection Stabilization
}
\author{Javier de Frutos\thanks{Instituto de Investigaci\'on en Matem\'aticas (IMUVA),
Universidad de Valladolid, Spain. Research supported
by Spanish MINECO
under grants MTM2013-42538-P (MINECO, ES) and MTM2016-78995-P (AEI/FEDER, UE).  (frutos@mac.uva.es)} \and Bosco
Garc\'{\i}a-Archilla\thanks{Departamento de Matem\'atica Aplicada
II, Universidad de Sevilla, Sevilla, Spain. Research supported by
Spanish MINECO under grant MTM2015-65608-P (bosco@esi.us.es)}
\and
Volker John\thanks{Weierstrass Institute for Applied Analysis and Stochastics,
Leibniz Institute in Forschungsverbund Berlin e. V. (WIAS), Mohrenstr. 39, 10117 Berlin, Germany.
Freie Universit\"at Berlin,
Department of Mathematics and Computer Science,
Arnimallee 6, 14195 Berlin, Germany.}
  \and Julia Novo\thanks{Departamento de
Matem\'aticas, Universidad Aut\'onoma de Madrid, Spain.  Research supported
by Spanish MINECO
under grants MTM2013-42538-P (MINECO, ES) and MTM2016-78995-P (AEI/FEDER, UE) (julia.novo@uam.es)}}
\date{\today}
\begin{document}
\maketitle
\abstract{This paper studies non inf-sup stable finite element approximations to the evolutionary Navier--Stokes equations.
Several local projection stabilization (LPS) methods corresponding to different stabilization
terms are analyzed, thereby separately studying the effects of the different stabilization terms.
Error estimates are derived in which the constants in the error bounds are independent
of inverse powers of the viscosity.
For one of the methods, using velocity and pressure finite elements of degree $l$, it will
be proved that the velocity error in $L^\infty(0,T;L^2(\Omega))$ decays with rate $l+1/2$
in the case that $\nu\le h$, with $\nu$ being the dimensionless viscosity and $h$ the mesh
width.
In the analysis of another method, it was observed that the convective
term can be bounded in an optimal way with the LPS stabilization of the pressure gradient.
Numerical studies confirm the analytical results.
}

\section{Introduction}

Let $\Omega \subset {\mathbb R}^d$, $d \in \{2,3\}$, be a bounded domain with
polyhedral  and Lipschitz boundary $\partial \Omega$. The incompressible
Navier--Stokes equations model the conservation of linear momentum and the
conservation of mass (continuity equation) by
\begin{align}
\label{NS} \partial_t\bu -\nu \Delta \bu + (\bu\cdot\nabla)\bu + \nabla p &= \bff &&\text{in }\ (0,T]\times\Omega,\nonumber\\
\nabla \cdot \bu &=0&&\text{in }\ (0,T]\times\Omega,\\
\bu(0, \cdot) &= \bu_0(\cdot)&&\text{in }\ \Omega,\nonumber
\end{align}
where $\bu$ is the velocity field, $p$ the kinematic pressure, $\nu>0$ the kinematic viscosity coefficient,
$\bu_0$ a given initial velocity, and $\bff$ represents the external body accelerations acting
on the fluid. The Navier--Stokes equations \eqref{NS} are equipped with homogeneous
Dirichlet boundary conditions $\bu = \boldsymbol 0$ on $\partial \Omega$.

This paper studies approximations to the Navier--Stokes equations (\ref{NS}) with non inf-sup stable
mixed finite elements in space and the implicit Euler method in time. We use the so-called local projection stabilization
(LPS) method to stabilize the pressure (since non inf-sup stable elements are used) plus other
stabilization terms which aim at allowing to derive
error estimates
where the constants do not depend explicitly on inverse powers of the viscosity but only implicitly through norms
of the solution of \eqref{NS}. This kind of bounds are called semi-robust or quasi-robust in the literature, see for example \cite{ADL16}.

In the literature, one can find already investigations of LPS methods for approximating
the solution of \eqref{NS}.
LPS methods for inf-sup stable elements are analyzed in  \cite{Lube_etalNS}.
The derived error bounds depend explicitly on inverse powers of the viscosity parameter $\nu$, unless the grids are
becoming sufficiently fine ($h \lesssim \sqrt\nu$, where $h$ is the mesh width), see also \cite{DA16} where error bounds for the {O}berbeck-{B}oussinesq model
model are obtained with an assumption on the regularity of the finite element solution.
In \cite{chacon_etal}, the authors consider non inf-sup stable mixed finite elements with LPS
stabilization. The so called term-by-term stabilization is applied, see \cite{Cha_maca}.  This method is a particular
type of a LPS method that is based on continuous functions, it does not need enriched finite element spaces,
and an interpolation operator replaces the standard projection operator
of the classical LPS methods. As in the present paper, a fully discrete scheme with the
implicit Euler method as time integrator is considered. A fully discrete LPS method for inf-sup stable
pairs of finite element spaces and a pressure-projection scheme
is analyzed in \cite{ADL16}.

Our analysis starts as in \cite{chacon_etal}, but there are   several major
differences in the formulation of the discrete problem as well as in the obtained results.
First of all, as an important result which was not achieved in
\cite{chacon_etal}, we are able to derive error bounds in which the constants do not depend on inverse powers of the diffusion parameter.
Also, contrary to \cite{chacon_etal}, where only one method is analyzed (with LPS stabilizations of the pressure,
the divergence, and the convective term),
we consider several methods, because  our aim is to study separately the effects of the
different stabilization terms. For all of them,  error bounds with constants independent on inverse powers of the
diffusion parameter are achieved
with the smallest possible number of stabilization terms. Also, in contrast to \cite{chacon_etal},
only moderate assumptions on the smallness of the time step $\Delta t$ are needed, like
$\Delta t\le C h^{d/2}$ in the error analysis of the pressure,
while in  \cite{chacon_etal} the smallness assumption on the mesh width $C h\le \Delta t$ is required.

Section~\ref{sec:global_grad_div} considers a method with LPS stabilization for the pressure and
a global grad-div stabilization term. The global grad-div stabilization term was proposed
to reduce the violation of mass conservation of finite element methods, but there are already
investigations which show that this term also stabilizes dominant convection.
In \cite{loscuatro_NS}, semi-robust error estimates are proved for
the standard Galerkin method plus grad-div stabilization in the case of inf-sup stable elements,
both for the continuous-in-time case and for the fully discrete case.
Paper \cite{loscuatro_NS} considers both, the regular case and the situation in which nonlocal compatibility
conditions for the solution are not assumed.
The results of Section~\ref{sec:global_grad_div} can be seen as an extension of some of the results
from \cite{loscuatro_NS} to the case of non inf-sup stable elements and also as an improvement of the
results from \cite{chacon_etal}. Error bounds of order $\mathcal O(h^s)$ are obtained for a sufficiently
smooth solution, where $2\le s\le l$, $s$ being the regularity index of the solution and $l$ being the degree of the polynomials used. The error is bounded in a
norm that includes the $L^2$ norm of the velocity at the final time step and the $L^2$ norm of the divergence.
This rate of convergence is the same as obtained in \cite{loscuatro_NS} for a similar norm and also
the same rate as proved in \cite{chacon_etal}. However, as we pointed out above, in  \cite{chacon_etal}
more terms are included in the method, the bound depends explicitly on $\nu^{-1}$, and the restriction
$C h\le \Delta t$ is assumed. For the error bound of the pressure, we get the optimal order
$\mathcal O(h^s)$. However, following the ideas of \cite{chacon_etal}, we are able to bound the
error of the $L^2$ norm of a discrete in time primitive of the pressure instead of the stronger
discrete in time $L^2$ norm of the pressure. Although Section~\ref{sec:global_grad_div} studies
the term-by-term stabilization, the analysis also holds for the standard one-level LPS method,
see \cite{Ga-Mat-To,Ma-Skr-To},  with slight modifications.

In Section~\ref{sec:lps_control_grad}, we analyze a method with LPS stabilization for the pressure
and LPS stabilization with control of the fluctuations of the gradient. For this section,
the use of term-by-term stabilization is necessary since in the error analysis
we need to have the same polynomial spaces for the velocity and the pressure. A key ingredient in the error
analysis is the application of \cite[Theorem~2.2]{berto}.
This result was already applied in the error analysis in \cite{Burman_Fer_numer_math}, where the authors
proved semi-robust error bounds for the evolutionary Navier--Stokes equations and a
continuous interior penalty (CIP) method in space assuming enough regularity of the solution. For the
method studied in Section~\ref{sec:lps_control_grad}, the convective term is estimated in an optimal
way (with constants independent on inverse powers of the diffusion parameter) with the help of the LPS stabilization of the gradient of the pressure. This LPS term was introduced
in \cite{BB01} to account for the violation of the discrete inf-sup condition by the used pair of finite elements.

Following the analysis of the previous section, Section~\ref{sec:lps_grad_div} presents
analogous error bounds for a method with both LPS stabilization for the pressure and the divergence.

For the methods analyzed in Sections~\ref{sec:global_grad_div} -- \ref{sec:lps_grad_div}, error
estimates with constants independent on inverse powers of the diffusion parameter
are derived with the help of stabilization terms that were not proposed for
stabilizing dominant convection but  to account for the non-satisfaction of
the discrete inf-sup condition or the violation of the
mass conservation (note that the LPS term of the velocity gradient of the method from Section~\ref{sec:lps_control_grad}
was not utilized for estimating the convective term).
The deeper reasons for this behavior are not yet understood and their explanation
is formulated as an open problem in \cite{paper-open}.

In Section~\ref{sec:ordersplusonehalf}, it is shown that the rate of decay of
the velocity
error in the situation $\nu \le h$ can be improved for the method from Section~\ref{sec:lps_control_grad} by
choosing different values of the stabilization parameters and increasing the regularity assumption
for the pressure. Concretely, a bound of order  $\mathcal O(h^{s+1/2})$ is proved for an error
which contains the $L^2$ error of the velocity. This is the same order that was obtained for the
CIP method in \cite{Burman_Fer_numer_math} under the same regularity assumptions. We are not aware
of any other paper where this order is proved and it is still an open question whether the
optimal expected order $\mathcal O(h^{s+1})$ for the $L^2$ error of the velocity can be achieved or not,
see \cite{paper-open}.

Finally, Section~\ref{sec:numres} presents numerical studies that confirm the analytical results.

\section{Preliminaries and notation}\label{sect:prelim}

Throughout the paper, $W^{s,p}(D)$ will denote the Sobolev space of real-valued functions defined on the domain $D\subset\mathbb{R}^d$ with distributional derivatives of order up to $s$ in $L^p(D)$. These spaces are endowed with the usual norm denoted by $\|\cdot\|_{W^{s,p}(D)}$.
 If $s$ is not a positive integer, $W^{s,p}(D)$ is defined by interpolation \cite{Adams}.
 In the case $s=0$ it is $W^{0,p}(D)=L^p(D)$. As it is standard, $W^{s,p}(D)^d$ will be endowed with the product norm and, since no confusion can arise, it will be denoted again by $\|\cdot\|_{W^{s,p}(D)}$. The case  $p=2$  will be distinguished by using $H^s(D)$ to denote the space $W^{s,2}(D)$. The space $H_0^1(D)$ is the closure in $H^1(D)$ of the set of infinitely differentiable functions with compact support
in $D$.  For simplicity, $\|\cdot\|_s$ (resp. $|\cdot |_s$) is used to denote the norm (resp. semi norm) both in $H^s(\Omega)$ or $H^s(\Omega)^d$. The exact meaning will be clear by the context. The inner product of $L^2(\Omega)$ or $L^2(\Omega)^d$ will be denoted by $(\cdot,\cdot)$ and the corresponding norm by $\|\cdot\|_0$.
For vector-valued functions, the same conventions will be used as before.
The norm of the dual space  $H^{-1}(\Omega)$  of $H^1_0(\Omega)$
is denoted by $\|\cdot\|_{-1}$.
As usual, $L^2(\Omega)$ is always identified
with its dual, so one has $H^1_0(\Omega)\subset L^2(\Omega)\subset H^{-1}(\Omega) $ with compact injection.
The following Sobolev's embedding \cite{Adams} will be used in the analysis: For  $1\le p<d/s$
let $q$ be such that $\frac{1}{q}
= \frac{1}{p}-\frac{s}{d}$. There exists a positive constant $C$, independent of $s$, such
that
\begin{equation}\label{sob1}
\|v\|_{L^{q'}(\Omega)} \le C \| v\|_{W^{s,p}(\Omega)}, \qquad
\frac{1}{q'}
\ge \frac{1}{q}, \quad v \in
W^{s,p}(\Omega).
\end{equation}
If $p>d/s$ the above relation is valid for $q'=\infty$. A similar embedding inequality holds for vector-valued functions.

Using the function spaces
$$
V=H_0^1(\Omega)^d,\quad Q=L_0^2(\Omega)=\left\{q\in L^2(\Omega):
(q,1)=0\right\},
$$
the weak formulation of problem (\ref{NS}) is as follows:
Find $(\bu,p)\in V\times Q$ such that for all $(\bv,q)\in V\times Q$,
\begin{equation}\label{eq:NSweak}
(\partial_t\bu,\bv)+\nu (\nabla \bu,\nabla \bv)+((\bu\cdot \nabla) \bu,\bv)-(\nabla \cdot \bv,p)
+(\nabla \cdot \bu,q)=(\boldsymbol f,\bv),
\end{equation}
and $\bu(0, \cdot)=\bu_0(\cdot)$.

The Hilbert space
$$
H^{\rm div}=\{ \bu \in L^{2}(\Omega)^d \ \mid \ L^2(\Omega) \ni \nabla \cdot \bu=0, \,
\bu\cdot \mathbf n|_{\partial \Omega} =0 \}$$
will be endowed with the inner product of $L^{2}(\Omega)^{d}$ and the space
$$V^{\rm div}=\{ \bu \in V \ \mid \ \nabla \cdot \bu=0 \}$$
with the inner product of $V$.

In the error analysis, the Poincar\'{e}--Friedrichs inequality
\begin{equation}\label{eq:poin}
\|\boldsymbol{v}\|_{0} \leq
C_{PF}\|\nabla \boldsymbol{v}\|_{0}\quad\forall v\in V
\end{equation}
will be used.

\section{Local projection stabilization with global grad-div stabilization.}\label{sec:global_grad_div}

Let $\mathcal T_h$  be a family of triangulations of $\overline\Omega$. Given an integer $l\ge 0$ and a mesh cell $K\in \mathcal T_h$ we
denote by $\Bbb P_l(K)$ the space of polynomials of degree less or equal to $l$. We consider the following finite element spaces
\begin{eqnarray*}
Y_h^l&=& \left\{v_h\in C^0(\overline\Omega)\mid {v_h}_{\mid_K}\in {\Bbb P}_l(K),\quad \forall K\in \mathcal T_h\right\}, \ l\ge 1,\nonumber\\
{\boldsymbol Y}_h^l&=&(Y_h^l)^d,\quad {\boldsymbol X}_h={\boldsymbol Y}_h^l\cap (H_0^1)^d, \nonumber\\
Q_h&=&Y_h^l\cap L_0^2.
\end{eqnarray*}
It will be assumed that the family of
meshes is quasi-uniform and that  the following inverse
inequality holds for each $v_{h} \in Y_h^l$, e.g., see \cite[Theorem 3.2.6]{Cia78},
\begin{equation}
\label{inv} \| \bv_{h} \|_{W^{m,p}(K)} \leq C_{\mathrm{inv}}
h_K^{n-m-d\left(\frac{1}{q}-\frac{1}{p}\right)}
\|\bv_{h}\|_{W^{n,q}(K)},
\end{equation}
where $0\leq n \leq m \leq 1$, $1\leq q \leq p \leq \infty$, and $h_K$
is the size (diameter) of the mesh cell $K \in \mathcal T_h$.

We consider the approximation of (\ref{NS}) with the implicit Euler method in time and a LPS method with grad-div stabilization in space.
Given ${\bu}_h^0=I_h \bu_0$, find $(\bu_h^{n+1},p_h^{n+1})\in {\boldsymbol X}_h\times Q_h$ such that
\begin{eqnarray}\label{eq:gal}
\lefteqn{\hspace*{-16em}\left(\frac{\bu_h^{n+1}-\bu_h^n}{\Delta t},\bv_h\right)+\nu(\nabla \bu_h^{n+1},\nabla \bv_h)+b(\bu_h^{n+1},\bu_h^{n+1},\bv_h)-(p_h^{n+1},\nabla \cdot \bv_h)}
\nonumber\\
+
S_h(\bu_h^{n+1},\bv_h)&=&({\boldsymbol f}^{n+1},\bv_h) \quad \forall \bv_h\in {\boldsymbol X}_h,\\
(\nabla \cdot \bu_h^{n+1},q_h)+s_{\rm pres}(p_h^{n+1},q_h)&=&0 \quad \forall q_h\in Q_h,\nonumber
\end{eqnarray}
where
\begin{eqnarray*}
S_h(\bu,\bv)&=&\mu(\nabla \cdot \bu,\nabla \cdot \bv),\\
b(\bu,\bv,\bw)&=&(B(\bu,\bv),\bw)\quad \forall\bu,\bv,\bw\in H_0^1(\Omega)^d,\\
B(\bu,\bv)&=&(\bu\cdot \nabla )\bv+\frac{1}{2}(\nabla  \cdot\bu)\bv\quad \forall \bu,\bv\in H_0^1(\Omega)^d,\\
s_{\rm pres}(p_h^{n+1},q_h)&=&\sum_{K\in \mathcal T_h}\tau_{p,K}(\sigma^*_h(\nabla p_h^{n+1}),\sigma^*_h(\nabla q_h))_K,
\end{eqnarray*}
and $\mu$ and $\tau_{p,K}$ are the grad-div and pressure stabilization parameters, respectively.
In addition,   $\sigma^*_h=Id-\sigma_h^{l-1}$, where $\sigma_h^{j}$ is a locally stable projection or interpolation operator from $L^2(\Omega)^d$ on ${\boldsymbol Y}_h^{j}$, that is,
there exists a constant $C>0$ such that for any
$K\in  {\mathcal T}_h$
\begin{eqnarray}\label{eq:stasigma}
\|\sigma_h^j(\bv)\|_{L^2(K)} \le C \|\bv\|_{L^2(\omega_K)},\quad \forall \bv\in {L^2(\Omega)}^d,
\end{eqnarray}
where $\omega_K$ is the union of all mesh cells whose intersection with $K$ is not empty. It will be assumed that the number of mesh cells in each set
$\omega_K$ is bounded independently of the triangulation and of $K$. From
\eqref{eq:stasigma}, also the $L^2$ stability of $\sigma^*_h$ follows.
The operator $\sigma_h^j$ can be chosen as a Bernardi--Girault~\cite{Ber_Gir}, Girault--Lions~\cite{Girault-Lions-2001}, or the Scott--Zhang~\cite{Scott-Z} interpolation operator
 in the space ${\boldsymbol Y}_h^{j}$ (for a proof of~(\ref{eq:stasigma}) in the case of the last two operators see~\cite{brenner-scot}). The following bound holds for $\bv\in H^s(\Omega)^d$,
 \begin{eqnarray}\label{eq:cota_sigma_local}
 \|\bv-\sigma_h^j(\bv)\|_{L^2(K)}\le C h_K^s\|\bv\|_{H^s(\omega_K)},\quad 1\le s\le j+1
 \end{eqnarray}
 from which it can be deduced that
\begin{eqnarray}\label{eq:cota_sigma}
\|\bv-\sigma_h^j(\bv)\|_0\le C h^{s}|\bv|_{s},\quad 1\le s\le j+1
\end{eqnarray}
 see \cite{Scott-Z,Ber_Gir,brenner-scot}. Bounds (\ref{eq:cota_sigma_local}) and (\ref{eq:cota_sigma}) will be applied for $j\in\left\{l-1,l\right\}$.


In the sequel, we will assume that
\begin{equation}\label{eq:tau_p}
\alpha_1h_K^2\le \tau_{p,K}\le \alpha_2 h_K^2
\end{equation}
for some positive constants $\alpha_1, \alpha_2$ independent of $h$.
In addition, the notations
\begin{equation}\label{eq:norm_tau_p}
(f,g)_{\tau_p}=\sum_{K\in {\mathcal T}_h}\tau_{p,K}(f,g)_K, \quad
\|f\|_{\tau_p} = (f,f)_{\tau_p}^{1/2}
\end{equation}
are used.

The following inf-sup condition holds (see \cite[Lemma 4.2]{chacon_etal}).
\begin{lema}\label{lema_inf_sup}
The following inf-sup condition holds
\begin{equation*}
\|q_h\|_0\le \beta_0\left(\sup_{\bv_h\in {\boldsymbol X}_h}\frac{(\nabla\cdot \bv_h,q_h)}{\|\nabla \bv_h\|_0}+
\|\sigma_h^*(\nabla q_h)\|_{\tau_p} \right) \quad \forall q_h \in Q_h.
\end{equation*}
\end{lema}
Along the paper we will use the following discrete Gronwall inequality whose proof can be found in \cite{hey-ran4}.
\begin{lema}\label{gronwall}
Let $k,B,a_j,b_j,c_j,\gamma_j$ be nonnegative numbers such that
$$
a_j+k\sum_{j=0}^n b_j\le k\sum_{j=0}^n \gamma_j a_j+k\sum_{j=0}^n c_j+B,\quad for \quad n\ge 0.
$$
Suppose that $k\gamma_j<1$, for all j, and set $\sigma_j=(1-k\gamma_j)^{-1}$. Then
$$
a_j+k\sum_{j=0}^nb_j\le \exp\left(k\sum_{j=0}^n\sigma_j\gamma_j\right)\left\{k\sum_{j=0}^nc_j+B\right\}, \quad for\quad  n\ge 0.
$$
\end{lema}

\subsection{Error bound for the velocity}
\label{se:velo_grad-div}
Let us denote by $\bu^n=\bu(\cdot,t_n)$ and by $p^n=p(\cdot,t_n)$.
Following \cite{chacon_etal}, we consider an approximation $\hat\bu_h^n=R_h \bu^n\in {\boldsymbol X}_h\subset {\boldsymbol Y}_h^l$ satisfying
\begin{equation}\label{eq:orto}
(\bu^n-\hat \bu_h^n,\bv_h)=0,\quad \forall \bv_h \in {\boldsymbol Y}_h^{l-1},\quad n=0,1,\ldots,N.
\end{equation}
Let us observe that the above definition for $\hat \bu_h$ can be applied for any time $t$ so that we can consider that $\hat \bu_h$ is continuous in the $t$ variable.
The following bound holds, see \cite{chacon_etal},
\begin{equation}\label{eq:hatu}
\|\bu^n-\hat \bu_h^n\|_{W^{m,p}}\le C h^{s+1-m+d/p-d/2}|\bu^n|_{s+1},\quad n=0,1,\ldots,N,
\end{equation}
for $m=0,1$, $p\in[1,\infty]$, $s\ge 1$.

Let $\hat p_h^n=I_h p^n\in Q_h$ with $I_h$ being the standard interpolation operator. There exists a constant $C>0$ such that
\begin{equation}\label{eq:hatp}
\|p^n-\hat p_h^n\|_{W^{m,p}}\le C h^{s-m+d/p-d/2}|p^n|_{s},\quad n=0,1,\ldots,N,\quad m=0,1,
\end{equation}
see \cite{brenner-scot}.

Let us denote
\begin{equation}\label{eq:err_split_vel}
\hat \be_h^{n}=\hat \bu_h^{n}-\bu^{n},\quad \be_h^n=\hat \bu_h^n-\bu_h^n,\quad \hat \lambda_h^{n}=\hat p_h^{n}-p^{n},\quad \lambda_h^n=\hat p_h^n-p_h^n.
\end{equation}
Subtracting the discrete problem \eqref{eq:gal} from the continuous problem \eqref{eq:NSweak}
yields the error equation
\begin{eqnarray}\label{eq:er}
\lefteqn{\left(\frac{\be_h^{n+1}-\be_h^n}{\Delta t},\bv_h\right)+\nu(\nabla \be_h^{n+1},\nabla \bv_h)+b(\hat \bu_h^{n+1},\hat\bu_h^{n+1},\bv_h)-b(\bu_h^{n+1},\bu_h^{n+1},\bv_h)}
\nonumber\\
&&-(\lambda_h^{n+1},\nabla \cdot \bv_h)+(\nabla \cdot \be_h^{n+1},q_h)+s_{\rm pres}(\lambda_h^{n+1},q_h)+
S_h(\be_h^{n+1},\bv_h)
\\
&=&({\boldsymbol \xi}_{v_h}^{n+1},\bv_h)+({\boldsymbol \xi}_{q_h}^{n+1},q_h)+
\nu(\nabla \hat \be_h^{n+1},\nabla \bv_h)
+s_{\rm pres}(\hat p_h^{n+1},q_h)\nonumber\\
&&+S_h(\hat \bu_h^{n+1},\bv_h)-(\hat \lambda_h^{n+1},\nabla \cdot \bv_h),\nonumber
\end{eqnarray}
for all $\bv_h\in {\boldsymbol X}_h$ and $q_h\in Q_h$. In \eqref{eq:er}, ${\boldsymbol \xi}_{v_h}^{n+1}$ and ${\boldsymbol \xi}_{p_h}^{n+1}$ are defined as follows
\begin{eqnarray}
\label{eq:trun1}
{\boldsymbol \xi}_{v_h}^{n+1}&=&{\boldsymbol \xi}_{v_h,1}^{n+1}+{\boldsymbol \xi}_{v_h,2}^{n+1},\\
\label{eq:trun11}
({\boldsymbol \xi}_{v_h,1}^{n+1},\bv_h)&=&-\left(\partial_t\bu^{n+1}-\frac{\hat \bu_h^{n+1}-\hat \bu_h^n}{\Delta t},\bv_h\right),
\\
({\boldsymbol \xi}_{v_h,2}^{n+1},\bv_h)&=&-b(\bu^{n+1},\bu^{n+1},\bv_h)+b(\hat \bu_h^{n+1},\hat \bu_h^{n+1},\bv_h),\label{eq:trun12}
\\
({\boldsymbol \xi}_{q_h}^{n+1},q_h)&=&(\nabla \cdot \hat \be_h^{n+1},q_h).
\nonumber 
\end{eqnarray}

\begin{remark}
Note that  the error equation (\ref{eq:er}) holds even for $(\bv_h,q_h)=(0,q_h)$ with $q_h \in Y_h^l$.  Let $q_h \in Y_h^l$ and denote by
$m(q_h)$ the mean of $q_h$, then  \eqref{eq:er} gives
\begin{eqnarray*}
\lefteqn{(\nabla \cdot \be_h^{n+1},q_h-m(q_h))+s_{\rm pres}(\lambda_h^{n+1},q_h-m(q_h))}
\\
&=&(\nabla \cdot \hat\be_h^{n+1},q_h-m(q_h))+s_{\rm pres}(\hat p_h^{n+1},q_h-m(q_h)).
\end{eqnarray*}
Since the terms $(\nabla \cdot \be_h^{n+1},m(q_h))$,  $s_{\rm pres}(\lambda_h^{n+1},m(q_h))$,
$(\nabla \cdot \hat\be_h^{n+1},m(q_h))$, and $s_{\rm pres}(\hat p_h^{n+1},m(q_h))$
vanish, it follows that
\begin{equation*}
(\nabla \cdot \be_h^{n+1},q_h)+s_{\rm pres}(\lambda_h^{n+1},q_h)
=(\nabla \cdot \hat\be_h^{n+1},q_h)+s_{\rm pres}(\hat p_h^{n+1},q_h) \quad \forall\ q_h\in Y_h^l.
\end{equation*}
\end{remark}

Setting
$(\bv_h,q_h)=(\be_h^{n+1},\lambda_h^{n+1})$, rearranging terms, and using the Cauchy--Schwarz
inequality and Young's inequality gives
\begin{eqnarray}\label{eq:unique}
\lefteqn{\frac{\|\be_h^{n+1}\|_0^2}{2\Delta t}-\frac{\|\be_h^{n}\|_0^2}{2\Delta t}+\frac{\|\be_h^{n+1}-\be_h^n\|_0^2}{2\Delta t}
+\frac{\nu}{2}\|\nabla \be_h^{n+1}\|_0^2
+\|\sigma_h^*(\nabla \lambda_h^{n+1})\|_{\tau_p}^2}
\nonumber\\
&&+S_h( \be_h^{n+1}, \be_h^{n+1})\nonumber\\
&\le& \left|b(\bu_h^{n+1},\bu_h^{n+1},\be_h^{n+1})-b(\hat \bu_h^{n+1},\hat \bu_h^{n+1},\be_h^{n+1})\right|+\frac{\|{\boldsymbol \xi}_{v_h}^{n+1}\|_0^2}{2}+\frac{\|\be_h^{n+1}\|_0^2}{2}\\
&&+\left|({\boldsymbol \xi}_{q_h}^{n+1},\lambda_h^{n+1})\right|
+\frac{\nu}{2}\|\nabla\hat \be_h^{n+1}\|_0
+\left|s_{\rm  pres}(\hat p_h^{n+1},\lambda_h^{n+1})\right|\nonumber\\
&&
+\left|S_h(\hat\bu_h^{n+1},\be_h^{n+1})\right|
+\left|(\hat \lambda_h^{n+1},\nabla \cdot \be_h)\right|.\nonumber
\end{eqnarray}

Now, the terms on the right-hand sid of \eqref{eq:unique} will be bounded. We start with
the last two terms. Applying the Cauchy--Schwarz inequality, Young's inequality, and
\eqref{eq:hatu} yields
\begin{align}\label{eq:erdiv}
\left|S_h(\hat\bu_h^{n+1},\be_h^{n+1})\right|&=\mu \left|(\nabla\cdot\hat\bu_h^{n+1},\nabla\cdot\be_h^{n+1})\right|
\le \frac{\mu}{8}\|\nabla\cdot\be_h^{n+1}\|_0^2 + 2\mu\|\nabla\cdot\hat \be_h^{n+1}\|_0^2
\nonumber\\
&\le \frac{1}{8}S_h(\be_h^{n+1}\,\be_h^{n+1})+ C\mu h^{2s}\|\bu\|_{L^\infty(H^{s+1})}^2.
\end{align}
Similarly, we obtain
\begin{equation}\label{eq:cotalambda}
\left|(\hat \lambda_h^{n+1},\nabla \cdot \be_h) \right| \le \frac{\mu}{8}\|\nabla \cdot \be_h^{n+1}\|_0^2
+\frac{2}{\mu} \|\hat \lambda_h^{n+1}\|_0^2
\le \frac{1}{8}S_h(\be_h^{n+1},\be_h^{n+1})+\frac{C}{\mu} h^{2s}\|p\|^2_{L^\infty(H^{s})},
\end{equation}
where in the last inequality \eqref{eq:hatp} was applied.
The nonlinear term in~(\ref{eq:unique}) can be bounded as in \cite{loscuatro_NS}
using the skew-symmetric property of $b$
\begin{eqnarray}\label{eq:ernonli}
\lefteqn{|b(\bu_h^{n+1},\bu_h^{n+1},\be_h^{n+1})-b(\hat \bu_h^{n+1},\hat \bu_h^{n+1},\be_h^{n+1})|}\nonumber\\
&\le&
|b(\be_h^{n+1},\hat\bu_h^{n+1},\be_h^{n+1})| +|b(\bu_h^{n+1}, \be_h^{n+1},\be_h^{n+1})|\nonumber\\
&\le& \|\nabla \hat \bu_h^{n+1}\|_{L^\infty}\|\be_h^{n+1}\|_0^2+\frac{1}{2}\|\nabla \cdot \be_h^{n+1}\|_0\|\hat \bu_h^{n+1}\|_{L^\infty}\|\be_h^{n+1}\|_0\nonumber\\
&\le& \left(\|\nabla \hat \bu_h^{n+1}\|_{L^\infty}+\frac{\| \hat \bu_h^{n+1}\|_{L^\infty}^2}{4\mu}\right)\|\be_h^{n+1}\|_0^2+\frac{\mu}{4}\|\nabla \cdot \be_h^{n+1}\|_0^2.
\end{eqnarray}
For the fourth term on the right-hand side of (\ref{eq:unique}), integrating by parts and using (\ref{eq:orto}), (\ref{eq:tau_p}), and (\ref{eq:hatu}) gives
\begin{eqnarray}\label{eq:erdiv2}
|({\boldsymbol \xi}_{q_h}^{n+1},\lambda_h^{n+1})|&=&|(\hat \be_h^{n+1},\nabla \lambda_h^{n+1})|
=|(\hat \be_h^{n+1},\sigma_h^*(\nabla \lambda_h^{n+1}))|
\nonumber\\&\le& \sum_{K\in {\mathcal T}_h}\frac{\|\hat  \be_h^{n+1}\|_{L^2(K)}^2}{\tau_{p,K}}+\frac{1}{4}\|\sigma_h^*(\nabla \lambda_h^{n+1}))\|_{\tau_p}^2\nonumber\\
&\le&C h^{2s}\|\bu\|_{L^\infty(H^{s+1})}^2+\frac{1}{4}\|\sigma_h^*(\nabla \lambda_h^{n+1}))\|_{\tau_p}^2.
\end{eqnarray}
For the fifth term, we use (\ref{eq:hatu}) to get
\begin{eqnarray}\label{eq:erehat}
\frac{\nu}{2}\|\nabla \hat \be_h^{n+1}\|_0^2\le C\nu h^{2s}\|\bu\|_{L^\infty(H^{s+1})}^2.
\end{eqnarray}
To bound the sixth  term, the usual inequalities, the definition \eqref{eq:norm_tau_p} of $\|\cdot\|_{\tau_p}$, \eqref{eq:cota_sigma},
and \eqref{eq:hatp} are utilized
\begin{eqnarray}\label{eq:erphat}
\left|s_{\rm pres}(\hat p_h^{n+1},\lambda_h^{n+1})\right|& \le& \|\sigma^*_h(\nabla \hat p_h^{n+1})\|^2_{\tau_{p}}+\frac{1}{4}\|\sigma^*_h(\nabla \lambda_h^{n+1})\|_{\tau_{p}}^2\nonumber\\
&\le& 2\|\sigma^*_h(\nabla \hat \lambda_h^{n+1})\|^2_{\tau_{p}}+2\|\sigma^*_h(\nabla p^{n+1})\|^2_{\tau_{p}}+\frac{1}{4}\|\sigma^*_h(\nabla \lambda_h^{n+1})\|_{\tau_{p}}\nonumber\\
&\le &  Ch^2\|\nabla \hat \lambda_h^{n+1}\|_0^2+C h^2\|\sigma^*_h(\nabla p^{n+1})\|_0^2+\frac{1}{4}\|\sigma^*_h(\nabla \lambda_h^{n+1})\|_{\tau_{p}}\nonumber\\
&\le& C h^{2s}\|p\|_{L^\infty(H^s)}^2+\frac{1}{4}\|\sigma^*_h(\nabla \lambda_h^{n+1})\|_{\tau_{p}}.
\end{eqnarray}
Inserting now (\ref{eq:erdiv}) -- (\ref{eq:erphat})  in (\ref{eq:unique}) yields
\begin{eqnarray*}
\lefteqn{\hspace*{-3em}\|\be_h^{n+1}\|_0^2-\|\be_h^n\|_0^2+\Delta t \nu\|\nabla \be_h^{n+1}\|_0^2+\Delta t\|\sigma_h^*(\nabla \lambda_h^{n+1})\|_{\tau_p}^2
+\mu\Delta t\|\nabla \cdot \be_h^{n+1}\|_0^2}\\
&\le&\Delta t \left(1+2\|\nabla \hat \bu_h^{n+1}\|_{L^\infty}+\frac{\|\hat\bu_h^{n+1}\|_{L^\infty}^2}{2\mu}\right)\|\be_h^{n+1}\|_0^2
+\Delta t\|{\boldsymbol\xi}_{v_h}^{n+1}\|_0^2\nonumber\\
&&+C \Delta t h^{2s}\left((1+\nu+\mu)\|\bu\|_{L^\infty(H^{s+1})}^2+\left(1+\mu^{-1}\right)\|p\|_{L^\infty(H^s)}^2\right),\nonumber
\end{eqnarray*}
such that summing over the discrete times leads to
\begin{eqnarray*}
\lefteqn{\|\be_h^{n}\|_0^2+\Delta t\nu\sum_{j=1}^n\|\nabla \be_h^j\|_0^2+\Delta t\sum_{j=1}^n\|\sigma_h^*(\nabla \lambda_h^j)\|_{\tau_p}^2
+\Delta t{\mu}\sum_{j=1}^n\|\nabla \cdot \be_h^j\|_0^2}\nonumber\\
&\le& \|\be_h^0\|_0^2+\sum_{j=1}^n\Delta t \left(1+2\|\nabla \hat \bu_h^{j}\|_{L^\infty}+\frac{\|\hat\bu_h^{j}\|_{L^\infty}^2}{2\mu}\right)\|\be_h^{j}\|_0^2+\Delta t\sum_{j=1}^n\|{\boldsymbol\xi}_{v_h}^j\|_0^2\nonumber\\
&&+C T h^{2s}\left((1+\nu+\mu)\|\bu\|_{L^\infty(H^{s+1})}^2+(1+\mu^{-1})\|p\|_{L^\infty(H^s)}^2\right).
\end{eqnarray*}

Let us bound $\|\hat\bu_h^{j}\|_{L^\infty}$ and $\|\nabla \hat \bu_h^{j}\|_{L^\infty}$, $1\le j\le n$.
For the first term, applying (\ref{sob1}) and (\ref{eq:hatu}) we have
\begin{eqnarray}\label{eq:cotahatu}
\|\hat\bu_h^j\|_{L^\infty}&\le& \|\bu^j\|_{L^\infty}+\|\bu^j-\hat \bu_h^j\|_{L^\infty}\le C\|\bu^j\|_{2}+Ch^{2-d/2}\|\bu^j\|_{2}\nonumber\\
&\le& C\|\bu\|_{L^\infty(H^{2})}.
\end{eqnarray}
Using the same argument for the second term, we reach
\begin{eqnarray}\label{eq:cotanablahatu}
\|\nabla \hat\bu_h^j\|_{L^\infty}&\le& \|\nabla \bu^j\|_{L^\infty}+\|\nabla \bu^j-\nabla \hat \bu_h^j\|_{L^\infty}\le C\|\bu^j\|_{3}+Ch^{2-d/2}\|\bu^j\|_{3}\nonumber\\
&\le& C\|\bu\|_{L^\infty(H^{3})}.
\end{eqnarray}
From (\ref{eq:cotahatu}) and (\ref{eq:cotanablahatu}) we deduce
\begin{eqnarray}\label{eq:cota_nonli}
1+2\|\nabla \hat \bu_h^{j}\|_{L^\infty}+\frac{\|\hat\bu_h^{j}\|_{L^\infty}^2}{2\mu}\le \hat M_u, \quad
\hat M_u=1+C\left(2\|\bu\|_{L^\infty(H^{3})}+\frac{\|\bu\|_{L^\infty(H^{2})}^2}{2\mu}\right).
\end{eqnarray}
Let us assume
\begin{equation}\label{eq:assume_Delta_t}
\Delta t \hat M_u\le \frac{1}{2}.
\end{equation}
Applying the Gronwall lemma, Lemma~\ref{gronwall}, we get
\begin{eqnarray}\label{eq:after_gron}
\lefteqn{
\|\be_h^{n}\|_0^2+\Delta t\nu\sum_{j=1}^n\|\nabla \be_h^j\|_0^2+\Delta t\sum_{j=1}^n\|\sigma_h^*(\nabla \lambda_h^j)\|_{\tau_p}^2
+\Delta t{\mu}\sum_{j=1}^n\|\nabla \cdot \be_h^j\|_0^2}\nonumber\\
&\le& e^{2T\hat M_u}\left( \|\be_h^0\|_0^2+\Delta t\sum_{j=1}^n\|{\boldsymbol\xi}_{v_h}^j\|_0^2\right)\\
&&+C e^{2T\hat M_u}\left(T h^{2s}\left((1+\nu+\mu)\|\bu\|_{L^\infty(H^{s+1})}^2+(1+\mu^{-1})\|p\|_{L^\infty(H^s)}^2\right)\right).\nonumber
\end{eqnarray}

To conclude the bound we are left with the task of getting a bound for the second term on the right-hand-side of (\ref{eq:after_gron}).
For the first term in the truncation error we write
\begin{eqnarray}\label{eq:trunca1_1}
\partial_t\bu^{j}-\frac{\hat\bu_h^{j}-\hat \bu^{j-1}_h}{\Delta t}&=&\left(\partial_t\bu^{j}-\frac{\bu^{j}- \bu^{j-1}}{\Delta t}\right)
+\left(\frac{\bu^{j}- \bu^{j-1}}{\Delta t}-\frac{\hat\bu_h^{j}-\hat \bu^{j-1}_h}{\Delta t}\right)\\
&=&\frac{1}{\Delta t}\int_{t_{j-1}}^{t_j}(t-t_{j-1})\partial_{tt}\bu(t)~dt+\frac{1}{\Delta t}\int_{t_{j-1}}^{t_j}\partial_t(\bu-\hat \bu_h)(t)~dt.\nonumber
\end{eqnarray}
Applying (\ref{eq:hatu}) and the Cauchy-Schwarz inequality, we reach
\begin{eqnarray}\label{eq:trunca1_2}
\left\|\partial_t\bu^{j}-\frac{\hat\bu_h^{j}-\hat \bu^{j-1}_h}{\Delta t}\right\|_0^2
\le C \Delta t \int_{t_{j-1}}^{t_j}\|\partial_{tt}\bu\|_{0}^2~dt+\frac{ h^{2s}}{\Delta t}\int_{t_{j-1}}^{t_j}\|\partial_t\bu(t)\|_{s}^2~dt.
\end{eqnarray}
 For the second term in the truncation error (\ref{eq:trun12}), we apply \cite[Lemma~2]{loscuatro_NS} to get
 \begin{eqnarray}\label{eq:trunca2_1}
\lefteqn{\sup_{\bphi\in L^2(\Omega)^d,\ \|\bphi\|_0=1}\left|b(\bu^{j},\bu^{j},\bphi)-b(\hat\bu_h^{j},\hat \bu_h^{j},\bphi)\right|}
\nonumber\\
&\le& C\left(\|\hat \bu_h^j\|_{L^\infty}+\|\nabla \cdot \hat \bu^j_h\|_{L^{2d/(d-1)}}+\|\bu^j\|_{2}\right)\|\bu^j-\hat \bu_h^j\|_{1}.
 \end{eqnarray}
To bound $\|\nabla \cdot \hat \bu^j_h\|_{L^{2d/(d-1)}}$ we use (\ref{sob1}) and (\ref{eq:hatu})
\begin{eqnarray}\label{eq:cotahatuL2d}
\|\nabla \cdot \hat \bu^j_h\|_{L^{2d/(d-1)}}&\le& \|\nabla \cdot \hat \bu^j\|_{L^{2d/(d-1)}}
+\|\nabla \cdot (\hat \bu^j_h-\bu^j)\|_{L^{2d/(d-1)}}\nonumber\\
&\le& C\|\bu^j\|_{2}+C h^{1/2}\|\bu^j\|_{2}\nonumber\\
& \le & C \|\bu\|_{L^\infty(H^{2})}.
\end{eqnarray}
Inserting (\ref{eq:cotahatu}) and (\ref{eq:cotahatuL2d}) in (\ref{eq:trunca2_1}) gives
  \begin{eqnarray}\label{eq:trunca2_2}
   \sup_{\bphi\in L^2,\ \|\bphi\|_0=1}\left| b(\bu^{j},\bu^{j},\bphi)-b(\hat\bu_h^{j},\hat \bu_h^{j},\bphi)\right|
 \le C\|\bu\|_{L^\infty(H^{2})}\|\bu^j-\hat \bu_h^j\|_{1}.
  \end{eqnarray}
 Then from (\ref{eq:trunca1_2}), (\ref{eq:trunca2_2}), and \eqref{eq:hatu} we get
 \begin{equation*}
 \Delta t\sum_{j=1}^n\|{\boldsymbol\xi}_{v_h}^j\|_0^2\le C T h^{2s}\left(\|\bu\|_{L^\infty(H^{2})}^2\|\bu\|_{L^\infty(H^{s+1})}^2+\|\partial_t\bu\|_{L^\infty(H^{s})}^2\right)+C (\Delta t)^2\int_{t_0}^{t_n}\|\partial_{tt}\bu\|_{0}^2~dt.
 \end{equation*}
 Inserting this inequality in (\ref{eq:after_gron}) and applying the triangle
 inequality to the splitting of the error \eqref{eq:err_split_vel}
 finishes the proof of the error
 estimate for the velocity.

\begin{Theorem}\label{thm:velo_bound}
Let the solution of \eqref{eq:NSweak} be sufficiently smooth in space and time,
such that all norms appearing in the formulation of this theorem are well defined,
and let the
time step be sufficiently small such that \eqref{eq:assume_Delta_t} holds. Then, the
following error bound holds for $2\le s\le l$:
 \begin{eqnarray}\label{eq:after_gron_2}
\lefteqn{\|\bu^n-\bu_h^{n}\|_0^2+\Delta t\nu\sum_{j=1}^n\|\nabla(\bu^j - \bu_h^j)\|_0^2+\Delta t\sum_{j=1}^n\|\sigma_h^*(\nabla (p^j- p_h^j))\|_{\tau_p}^2
}\nonumber\\
&& +\Delta t{\mu}\sum_{j=1}^n\|\nabla \cdot \bu_h^j\|_0^2 \\
&\le& C e^{2T\hat M_u}\left( \|\be_h^0\|_0^2+  T \hat K_{u,p} h^{2s}+ (\Delta t)^2\int_{t_0}^{t_n}\|\partial_{tt}\bu\|_{0}^2~dt\right),\nonumber
\end{eqnarray}
where $\hat M_u$ is defined in \eqref{eq:cota_nonli} and
$$
\hat K_{u,p}=\left((1+\|\bu\|_{L^\infty(H^{2})}^2+\nu+\mu)\|\bu\|_{L^\infty(H^{s+1})}^2+\|\partial_t\bu\|_{L^\infty(H^{s})}^2+(1+\mu^{-1})\|p\|_{L^\infty(H^s)}^2\right).
$$
\end{Theorem}

Note that neither $\hat M_u$ nor $\hat K_{u,p}$ depend explicitly on negative powers of $\nu$.
The error bound \eqref{eq:after_gron_2} can be summarized in the form
$$
\mbox{errors on the left-hand side of \eqref{eq:after_gron_2} }\le
C(\bu,\partial_t\bu,\partial_{tt}\bu, p,T,\mu,\mu^{-1}) \left(\|\be_h^0\|_0+h^s+\Delta t\right).
$$

\subsection{Error bound for the pressure}
\label{sec:pres_bound_0}

We will derive now a bound for the error in the pressure. Let us denote
$$
\Lambda_h^n=\Delta t\sum_{j=1}^n\lambda_h^j,\quad \hat\Lambda_h^n=\Delta t\sum_{j=1}^n\hat\lambda_h^j.
$$
Setting $q_h=0$ in the error equation \eqref{eq:er} yields
\begin{eqnarray}\label{eq:Lambda}
(\Lambda_h^n,\nabla \cdot \bv_h) &=& (\be_h^{n}-\be_h^0,\bv_h)+\Delta t\nu\sum_{j=1}^n(\nabla (\bu^{j}-\bu_h^{j}),\nabla\bv_h)
\nonumber\\
&&+\Delta t\sum_{j=1}^n\left(b(\bu^{j},\bu^{j},\bv_h)-b(\bu_h^{j},\bu_h^{j},\bv_h)\right) +\Delta t\mu\sum_{j=1}^n(\nabla \cdot(\bu^{j}-\bu_h^{j}),\nabla \cdot \bv_h)\nonumber \\
&& +(\hat \Lambda_h^n,\nabla \cdot \bv_h)
+\Delta t\sum_{j=1}^n\left(\partial_t\bu^{j}-\frac{\hat \bu_h^{j}-\hat \bu_h^{j-1}}{\Delta t},\bv_h\right).
\end{eqnarray}
Applying Lemma~\ref{lema_inf_sup} we obtain
\begin{eqnarray}\label{eq:infsup_Lambda}
\|\Lambda_h^n\|_0\le \beta_0\left(\sup_{\bv_h\in {\boldsymbol X}_h}\frac{(\Lambda_h^n,\nabla\cdot \bv_h)}{\|\nabla \bv_h\|_0}+
\|\sigma_h^*(\nabla \Lambda_h^n)\|_{\tau_p} \right).
\end{eqnarray}
Let us bound the first term on the right-hand side of \eqref{eq:infsup_Lambda}. From (\ref{eq:Lambda}) we get with the triangle inequality, the Poincar\'e--Friedrichs
inequality \eqref{eq:poin}, and the
estimate for the dual pairing
\begin{eqnarray}\label{eq:infsup_Lmabda_2}
\sup_{\bv_h\in {\boldsymbol X}_h}\frac{(\Lambda_h^n,\nabla\cdot \bv_h)}{\|\nabla \bv_h\|_0}
&\le& \|\be_h^n\|_{-1}+\|\be_h^0\|_{-1}+\Delta t\nu\sum_{j=1}^n\|\nabla (\bu^j-\bu_h^j)\|_0\nonumber\\
&&+\Delta t\sum_{j=1}^n\|B(\bu^j,\bu^j)-B(\bu_h^j,\bu_h^j)\|_{-1}+\Delta t \mu\sum_{j=1}^n\|\nabla \cdot \bu_h^j\|_0\nonumber\\
&& + \Delta t \sum_{j=1}^n\|\hat\lambda_h^j\|_0+\Delta t \sum_{j=1}^n \left\|\partial_t\bu^j-\frac{\hat \bu_h^j-\hat\bu_h^{j-1}}{\Delta t}\right\|_{-1}.
\end{eqnarray}
Note that, since $\|\cdot\|_{-1} \le C \|\cdot\|_0$,  the first term on the right-hand side of \eqref{eq:infsup_Lmabda_2} was already bounded in the derivation of the velocity error bound. To bound the third and fifth
term on the right-hand side of~(\ref{eq:infsup_Lmabda_2}), we use the
fact that for any sequence $\left\{\alpha_j\right\}_{j=1}^\infty$ of nonnegative real numbers and $n\le T/\Delta t$ by the Cauchy--Schwarz inequality holds
\begin{equation}\label{eq:cota_gen}
 \Delta t\sum_{j=1}^n\alpha_j \le T^{1/2}\left(\Delta t\sum_{j=1}^n \alpha_j^2\right)^{1/2}.
\end{equation}
With this estimate and the velocity error bound \eqref{eq:after_gron_2}, an estimate for the
third and fifth term is obtained.
Using (\ref{eq:cota_gen}) and~(\ref{eq:trunca1_2}), the bound of the last term on the right-hand side of~(\ref{eq:infsup_Lmabda_2}) follows.
 For the sixth term, we apply (\ref{eq:hatp}) to get
$$
\Delta t\sum_{j=1}^n\|\hat\lambda_h^j\|_0\le C T h^{s}\|p\|_{L^\infty(H^{s})}.
$$
We are left with the fourth term on the right-hand side of (\ref{eq:infsup_Lmabda_2}).
Arguing as in~\cite{loscuatro_NS}, we obtain
\begin{eqnarray*}
\lefteqn{\Delta t\sum_{j=1}^n\|B(\bu^j,\bu^j)-B(\bu_h^j,\bu_h^j)\|_{-1}}\nonumber\\
&\le&  C\Delta t\sum_{j=1}^n\left(\|\bu_h^j\|_{L^\infty}+\|\nabla \cdot  \bu_h^j\|_{L^{2d/(d-1)}}+\|\bu^j\|_{2}\right)\|\bu^j-\bu_h^j\|_0\nonumber\\
&&+C\Delta t\sum_{j=1}^n\|\bu^j\|_{1}\|\nabla \cdot (\bu^j-\bu_h^j)\|_0\nonumber\\
&\le& CT\left( \max_{1\le j\le n}(\|\bu_h^j\|_{L^\infty}+\|\bu^j\|_2\right)\max_{1\le j\le n}\|\bu^j-\bu_h^j\|_0\nonumber\\
&& +CT^{1/2}\left(\Delta t \sum_{j=1}^n\|\nabla \cdot  \bu_h^j\|_{L^{2d/(d-1)}}^2\right)^{1/2}
\max_{1\le j\le n}\|\bu^j-\bu_h^j\|_0\nonumber\\
&& + C T^{1/2}\|\bu\|_{L^\infty(H^1)} \left(\Delta t \sum_{j=1}^n\|\nabla \cdot (\bu^j-\bu_h^j)\|_0^2\right)^{1/2}.
\end{eqnarray*}
To bound the norms involving $\bu_h^j$, the inverse inequality \eqref{inv}, the Sobolev
embedding \eqref{sob1}, and \eqref{eq:cotahatu} are used to get
\begin{eqnarray}\label{eq:norma_inf_uh}
\|\bu_h^j\|_{L^\infty}&\le& \|\be_h^j\|_{L^\infty}+\|\hat \bu_h^j\|_{L^\infty}
\le Ch^{-d/2}\|\be_h^j\|_0+\|\hat \bu_h^j\|_{L^\infty} \nonumber\\
&\le& Ch^{-d/2}\|\be_h^j\|_0+\|\bu^j - \hat \bu_h^j\|_{L^\infty} + \|\bu^j \|_{L^\infty}\nonumber\\
&\le&  C h^{-d/2} \|\be_h^j\|_0 + C h^{2-d/2} \|\bu\|_{2} +  C \|\bu\|_{2} \nonumber\\
&\le& C h^{-d/2} \|\be_h^j\|_0 + C\|\bu\|_{L^\infty(H^{2})}.
\end{eqnarray}
The term $\|\be_h^j\|_0$ was already bounded during the derivation of the velocity error
estimate. Applying the inverse estimate \eqref{inv} gives
\begin{equation}\label{eq:2d_dmenos1}
\left(\Delta t \sum_{j=1}^n\|\nabla \cdot  \bu_h^j\|_{L^{2d/(d-1)}}^2\right)^{1/2 } \le Ch^{-1/2} \left(\Delta t\sum_{j=1}^n\|\nabla \cdot  \bu_h^j\|_0^2\right)^{1/2},
\end{equation}
where the term on the right-hand side is already bounded in \eqref{eq:after_gron_2}.
Using \eqref{eq:norma_inf_uh}, \eqref{eq:2d_dmenos1} and assuming
\begin{equation}\label{extra_deltat}
\|\be_h^0\|_0 = \mathcal O(h^{d/2}) \quad \mbox{and} \quad
\Delta t \le C h^{d/2},
\end{equation}
 we finally reach
\begin{eqnarray*}
\lefteqn{\Delta t\sum_{j=1}^n\|B(\bu^j,\bu^j)-B(\bu_h^j,\bu_h^j)\|_{-1}}\nonumber\\
&\le&
C(\bu,\partial_t\bu,\partial_{tt}\bu,p,T,\mu,\mu^{-1}) \left(\max_{1\le j\le n}\|\bu^j-\bu_h^j\|_0+\left(\Delta t\sum_{j=1}^n \|\nabla\cdot (\bu^j-\bu_h^j)\|_0^2\right)^{1/2}\right).
\end{eqnarray*}
The bound of this term is finished by applying  \eqref{eq:after_gron_2}.

Inserting the derived inequalities in (\ref{eq:infsup_Lmabda_2}) and going back to (\ref{eq:infsup_Lambda}) yields
\begin{eqnarray*}
\|\Lambda_h^n\|_0&\le& \beta_0C(\bu,\partial_t\bu,\partial_{tt}\bu, p,T,\mu,\mu^{-1}) \left(\|\be_h^0\|_0+h^s+\Delta t\right)+\beta_0
\|\sigma_h^*(\nabla \Lambda_h^n)\|_{\tau_p}.
\end{eqnarray*}
The last term was already bounded in the derivation of the velocity error estimate, since it
is by the Cauchy--Schwarz inequality
\begin{eqnarray*}
\|\sigma_h^*(\nabla \Lambda_h^n)\|_{\tau_p}^2&=& \left\|\Delta t\sum_{j=1}^n \sigma_h^*(\nabla \lambda_h^j)\right\|_{\tau_p}^2\le
n(\Delta t)^2 \sum_{j=1}^n\| \sigma_h^*(\nabla \lambda_h^j)\|_{\tau_p}^2 \\
&=&T\Delta t \sum_{j=1}^n \|\sigma_h^*(\nabla \lambda_h^j)\|_{\tau_p}^2,
\end{eqnarray*}
which is a term on the left-hand side of estimate  \eqref{eq:after_gron}.
The estimate for the pressure error is obtained by applying finally the triangle
inequality to the splitting $p^j-p_h^j = \lambda_h^j - \hat   \lambda_h^j$ and
using \eqref{eq:hatp}.

\begin{Theorem} \label{thm:pres_bound}
Let the assumption of Theorem~\ref{thm:velo_bound} and the assumptions
\eqref{extra_deltat} be satisfied, then the following error estimate holds
\begin{equation*}
\left\|\Delta t \sum_{j=1}^n (p^j-p_h^j)\right\|_0
\le \beta_0C(\bu,\partial_t\bu,\partial_{tt}\bu, p,T,\mu^{-1}) \left(\|\bu_0 - \bu_h^0\|_0+h^s+\Delta t\right).
\end{equation*}
\end{Theorem}

\section{Local projection stabilization with control of the fluctuation of the gradient}
\label{sec:lps_control_grad}

In this part we will concentrate on the LPS method based on the stabilization of the gradient.
The stabilization term $S_h$ is defined by
\begin{equation}\label{eq:lps_gradient}
S_h(\bu_h,\bv_h) := \sum_{K\in \mathcal{T}_h}
\tau_{\nu,K}\left(\sigma^*_h (\nabla \bu_h), \sigma_h^*(\nabla \bv_h)\right)_K,
\end{equation}
where $\tau_{\nu,K}$, $K\in\mathcal{T}_h$, are non-negative constants.
This kind of  LPS method gives additional control on the fluctuation of the gradient.
In the sequel we will use the notations
$$
(f,g)_{\tau_\nu}=\sum_{K\in {\mathcal T}_h}\tau_{\nu,K}(f,g)_K\quad \mbox{and}\quad
\|f\|_{\tau_\nu} = (f,g)_{\tau_\nu}^{1/2}.
$$
For the stabilization parameter we will take $\tau_{\nu,K}
 \sim 1$. The same finite element spaces are used as in Section~\ref{sec:global_grad_div}.

\paragraph{\bf Assumption~A1}
There exits an interpolation operator $i_h\ :\
H^2(\Omega)\rightarrow Q_h$ with the approximation properties
\begin{equation}
\big\|q-i_hq\big\|_{0,K} + h_K \big| q-i_hq \big|_{1,K} \le Ch_K^{s+1} \big\|q\big\|_{s+1,K}
\quad \forall\  q\in H^{s+1}\big(K\big),\,1\le s\le l,\label{j2}
\end{equation}
for all $K\in\mathcal{T}_h$.
The pressure interpolation operator $i_h$ satisfies the orthogonality condition
\begin{align}\label{j3}
 (q-i_hq, r_h)_K = 0 \quad \forall\  q\in Q\cap H^2(\Omega),\ r_h\in Y_h^{l-1}, K\in\mathcal{T}_h.
\end{align}

\begin{remark}
The operator $i_h$ is the analog in the pressure space to the approximation used in the previous section to bound the velocity error.
\end{remark}\smallskip

Let us observe that the velocity and pressure spaces ${\boldsymbol Y}_h$ and $Y_h$, respectively,
are based on piecewise polynomials of the same degree $l$ and are the same space (apart from the fact that the velocity space has $d$ components). This property is essential for applying the following lemma. This lemma  can be deduced from \cite[Theorem~2.2]{berto}.

\begin{lema}\label{cor_burman}
Let $\sigma_h^j\ : \ L^2(\Omega)^d\to {\boldsymbol Y}_h^{j}$ be the interpolation operation defined
in Section~\ref{sec:global_grad_div} and let $\bu\in W^{1,\infty}(\Omega)^d$ and $\bv_h\in {\boldsymbol Y}_h^{j}$. Then, it holds
\begin{eqnarray}
\|(I-\sigma_h^j)(\bu\cdot \bv_h)\|_0&\le& C h \|\bu\|_{W^{1,\infty}}\|\bv_h\|_0,\nonumber \\ 
\|(I-\sigma_h^j)(\bu\cdot \bv_h)\|_{1}&\le& C \|\bu\|_{W^{1,\infty}}\|\bv_h\|_0\label{eq:berto}.
\end{eqnarray}
\end{lema}
Lemma~\ref{cor_burman} will be applied for $j\in \left\{l-1,l\right\}$.


\begin{remark}
Lemma~\ref{cor_burman} holds true for $\bv_h\in {\boldsymbol Y}_h^{j}$ with several components or $v_h\in Y_h^j$ with only one component.
\end{remark}\smallskip

\begin{remark}
In this section, in order to apply Lemma~\ref{cor_burman},
we need that the velocity and pressure spaces are the same. Then, the analysis holds for the LPS method based on the term-by-term stabilization introduced in \cite{Cha_maca}.
On the contrary, the analysis of the previous section also holds for the standard one-level LPS method over triangular or quadrilateral elements \cite{Ga-Mat-To,Ma-Skr-To} with slight modifications.
\end{remark}

\subsection{Error bound for the velocity}
We consider the approximation of (\ref{eq:NSweak}) with the implicit Euler method in time and a LPS method with LPS stabilization for the gradient of the velocity \eqref{eq:lps_gradient} and for the pressure.
Given ${\bu}_h^0=I_h \bu_0$, find $(\bu_h^{n+1},p_h^{n+1})\in ({\boldsymbol X}_h,Q_h)$ such that
\begin{eqnarray*}
\lefteqn{\hspace*{-16em}\left(\frac{\bu_h^{n+1}-\bu_h^n}{\Delta t},\bv_h\right)+b(\bu_h^{n+1},\bu_h^{n+1},\bv_h)+\nu(\nabla \bu_h^{n+1},\nabla \bv_h)
-(p_h^{n+1},\nabla \cdot \bv_h)}
\nonumber\\
+
S_h(\bu_h^{n+1},\bv_h)&=&({\boldsymbol f}^{n+1},\bv_h),\quad \forall \bv_h\in {\boldsymbol X}_h,\\
(\nabla \cdot \bu_h^{n+1},q_h)+s_{\rm pres}(p_h^{n+1},q_h)&=&0,\quad \forall q_h\in Q_h.\nonumber
\end{eqnarray*}

In the sequel, we will denote by $\hat \bu_h^n$ the function defined in Assumption~A1 satisfying (\ref{eq:orto}) and by $\hat p_h^n=i_h p^n$ and we denote
$$
\hat \be_h^{n}=\hat \bu_h^{n}-\bu^{n},\quad \be_h^n=\hat \bu_h^n-\bu_h^n,\quad \hat \lambda_h^{n}=\hat p_h^{n}-p^{n},\quad \lambda_h^n=\hat p_h^n-p_h^n.
$$
It is easy to see that $(\be_h^n, \lambda_h^{n})$  satisfies the same equation~(\ref{eq:er}) as in
Section~\ref{se:velo_grad-div} and,
consequently, (\ref{eq:unique}). In the present analysis, the
first term on the right-hand
side of (\ref{eq:unique}) and the last three ones will be treated differently.

Starting as for deriving \eqref{eq:ernonli} yields
\begin{eqnarray}\label{eq:esta}
\lefteqn{|b(\bu_h^{n+1},\bu_h^{n+1},\be_h^{n+1})-b(\hat \bu_h^{n+1},\hat \bu_h^{n+1},\be_h^{n+1})|}
\nonumber\\
&\le& \|\nabla \hat \bu_h^{n+1}\|_{L^\infty}\|\be_h^{n+1}\|_0^2+\frac{1}{2}((\nabla \cdot \be_h^{n+1})\hat \bu_h^{n+1},\be_h^{n+1}).
\end{eqnarray}
To bound the second term on the right-hand side of (\ref{eq:esta}), we decompose
\begin{eqnarray}\label{eq:likelu1}
&&\lefteqn{((\nabla \cdot \be_h^{n+1})\hat \bu_h^{n+1},\be_h^{n+1})}\\
&&=\left((\nabla \cdot \be_h^{n+1}),{{\sigma_h^l}}(\hat \bu_h^{n+1}\cdot \be_h^{n+1})\right)+\left((\nabla \cdot \be_h^{n+1}),(I-{\sigma_h^l})(\hat \bu_h^{n+1}\cdot \be_h^{n+1})\right).\nonumber
\end{eqnarray}
Using the error equation {\eqref{eq:er}} with $(\bv_h,q_h) = \left(\boldsymbol 0,{\sigma_h^l}(\hat \bu_h^{n+1}\cdot \be_h^{n+1})\right)$ gives
for the first term on the right-hand side of (\ref{eq:likelu1})
\begin{eqnarray}\label{eq:52}
\lefteqn{\left((\nabla \cdot \be_h^{n+1}),{\sigma_h^l}(\hat \bu_h^{n+1}\cdot\be_h^{n+1})\right)}\nonumber\\
&=&s_{\rm pres}(p_h^{n+1},{\sigma_h^l}(\hat \bu_h^{n+1}\cdot\be_h^{n+1}))
+(\nabla \cdot \hat \be_h^{n+1},{\sigma_h^l}(\hat \bu_h^{n+1}\cdot\be_h^{n+1})).
\end{eqnarray}
For the first term on the right-hand side in \eqref{eq:52},
arguing as in (\ref{eq:erphat}),  we have
\begin{eqnarray*}
\lefteqn{s_{\rm pres}(p_h^{n+1},{\sigma_h^l}(\hat \bu_h^{n+1}\cdot\be_h^{n+1}))}\nonumber\\
&\le& C h^{2s}\|p\|_{L^\infty(H^s)}^2+\frac{1}{8}\|\sigma^*_h(\nabla \lambda_h^{n+1})\|_{\tau_p}^2+C h^2\|\sigma^*_h(\nabla {\sigma_h^l}(\hat \bu_h^{n+1}\cdot\be_h^{n+1}))\|_{0}^2.
\end{eqnarray*}
For the last term above, applying \eqref{j2}, the inverse estimate \eqref{inv}, and {(\ref{eq:stasigma})},  it follows that
\begin{eqnarray*}
h^2\|\sigma^*_h(\nabla {\sigma_h^l}(\hat \bu_h^{n+1}\cdot\be_h^{n+1}))\|_{0}^2&\le& Ch^2 \|\nabla {\sigma_h^l}(\hat \bu_h^{n+1}\cdot\be_h^{n+1})\|_{0}^2\\
&\le& Ch^2 h^{-2}\|{\sigma_h^l}(\hat \bu_h^{n+1}\cdot\be_h^{n+1})\|_{0}^2\nonumber\\
&\le& C \|\hat \bu_h^{n+1}\|_{L^\infty}^2\|\be_h^{n+1}\|_{0}^2,
\end{eqnarray*}
so that
\begin{eqnarray}\label{eq:handle00}
\lefteqn{s_{\rm pres}(p_h^{n+1},{\sigma_h^l}(\hat \bu_h^{n+1}\cdot\be_h^{n+1}))}\nonumber\\
&\le&C h^{2s}\|p\|_{L^\infty(H^s)}^2+\frac{1}{8}\|\sigma^*_h(\nabla \lambda_h^{n+1})\|_{\tau_p}^2+ C \|\hat \bu_h^{n+1}\|_{L^\infty}^2\|\be_h^{n+1}\|_{0}^2.
\end{eqnarray}
To bound the second term on the right-hand side of (\ref{eq:52}), we get with \eqref{eq:hatu}
and {(\ref{eq:stasigma})}
\begin{eqnarray*}
(\nabla \cdot \hat \be_h^{n+1},{\sigma_h^l}(\hat \bu_h^{n+1}\cdot\be_h^{n+1}))\le C h^{2s}\|\bu\|_{L^\infty(H^{s+1})}^2+C \|\hat \bu_h^{n+1}\|_{L^\infty}^2\|\be_h^{n+1}\|_{0}^2.
\end{eqnarray*}
For the second term on the right-hand side of (\ref{eq:likelu1}), we apply Lemma~\ref{cor_burman} and the inverse inequality (\ref{inv}) to obtain
\begin{eqnarray*}
\left((\nabla \cdot \be_h^{n+1}),(I-{\sigma_h^l})(\hat u_h^{n+1}\cdot \be_h^{n+1})\right)&\le& C h \|\nabla \cdot \be_h^{n+1}\|_0 \|\hat \bu_h^{n+1}\|_{W^{1,\infty}}\|\be_h^{n+1}\|_0\nonumber\\
&\le& C \|\hat \bu_h^{n+1}\|_{W^{1,\infty}} \|\be_h^{n+1}\|_0^2.
\end{eqnarray*}
Collecting all estimates, we reach
\begin{eqnarray}\label{eq:ernonli2}
\lefteqn{
\left|b(\bu_h^{n+1},\bu_h^{n+1},\be_h^{n+1})-b(\hat \bu_h^{n+1},\hat \bu_h^{n+1},\be_h^{n+1})\right|}\nonumber\\
&\le&
C\left(\|\nabla \hat \bu_h^{n+1}\|_{L^\infty}+\|\hat \bu_h^{n+1}\|_{L^\infty}^2 \right)\|\be_h^{n+1}\|_0^2
+C h^{2s}\left(\|p\|_{L^\infty(H^s)}^2+\|\bu\|_{L^\infty(H^{s+1})}^2\right)\nonumber\\
&&+\frac{1}{8}\|\sigma^*_h(\nabla \lambda_h^{n+1})\|_{\tau_p}^2.
\end{eqnarray}
\medskip

\begin{remark} We like to emphasize the aspect that the only stabilization that was used
to derive the optimal estimate \eqref{eq:ernonli2} of the convective term (in which the constants do not depend on inverse powers
of the diffusion parameter) was the LPS
stabilization of the pressure -- a stabilization term whose proposal does not possess any
connection with dominant convection.
\end{remark}
\medskip

The last three terms on the right-hand side of~(\ref{eq:unique}) will be bounded next. The term $s_{\rm  pres}(\hat p_h^{n+1},\lambda_h^{n+1})$ can be bounded as in~(\ref{eq:erphat}),
using \eqref{j2} instead of \eqref{eq:hatp}, and replacing the factor 1/4 multiplying the last
term in~(\ref{eq:erphat}) by~1/8. Also,
arguing similarly to~(\ref{eq:erdiv}) we
have
\begin{eqnarray*}
S_h(\hat\bu_h^{n+1},\be_h^{n+1})&= &(\sigma^*_h(\nabla\hat\bu_h^{n+1}),\sigma^*_h(\nabla\be_h^{n+1}))_{\tau_\nu}\\
&\le& \frac{1}{4}\|\sigma^*_h(\nabla\be_h^{n+1})\|_{\tau_\nu}^2 + \|\sigma^*_h(\nabla\hat\bu_h^{n+1})\|_{\tau_\nu}^2
\nonumber\\
&= &\frac{1}{4}S_h(\be_h^{n+1} ,\be_h^{n+1})+  \|\sigma^*_h(\nabla\hat\bu_h^{n+1})\|_{\tau_\nu}^2.
\end{eqnarray*}
Then, applying the $L^2$ stability of $\sigma_h^*$, \eqref{eq:cota_sigma}, and (\ref{eq:hatu}) yields
\begin{eqnarray*}
\|\sigma_h^*(\nabla \hat\bu_h^{n+1})\|_{\tau_{\nu}}^2&\le& \|\sigma_h^*\nabla (\hat\bu_h^{n+1}-\bu^{n+1})\|_{\tau_{\nu}}^2+
\|\sigma_h^*(\nabla \bu^{n+1})\|_{\tau_{\nu}}^2\nonumber\\
&\le& C\|\nabla (\hat\bu_h^{n+1}-\bu^{n+1})\|_0^2+C h^{2s}\|\bu\|_{L^\infty(H^{s+1})}^2\nonumber\\
&\le& C h^{2s}\|\bu\|_{L^\infty(H^{s+1})}^2,
\end{eqnarray*}
so that
\begin{equation}
\label{eq:erdiv_4.1}
S_h(\hat\bu_h^{n+1},\be_h^{n+1})\le \frac{1}{4}S_h(\be_h^{n+1},\be_h^{n+1})+  C h^{2s}\|\bu\|_{L^\infty(H^{s+1})}^2.
\end{equation}

Finally, to bound the last term on the right-hand of (\ref{eq:unique}),  we use the orthogonality condition of the pressure interpolation operator (\ref{j3}), that the norm of the gradient contains all terms
of the norm of the divergence and
$\|\sigma_h^*(\nabla \cdot \be_h^{n+1})\|_{\tau_{\nu}} \le \sqrt{d}\|\sigma_h^*(\nabla \be_h^{n+1})\|_{\tau_{\nu}}$ holds, that
$\tau_{\nu,K}\sim 1$, and \eqref{j2} to get
\begin{eqnarray}\label{eq:last_term_4.1}
(\hat \lambda_h^{n+1},\nabla \cdot \be_h^{n+1})&=&-(p^{n+1}-i_h p^{n+1},\nabla\cdot \be_h^{n+1})=
-(p^{n+1}-i_h p^{n+1}, \sigma_h^*(\nabla\cdot  \be_h^{n+1}))\nonumber\\
&\le&\|p^{n+1}-i_h p^{n+1}\|_{\tau_{\nu}^{-1}}\|\sigma_h^*(\nabla \cdot \be_h^{n+1})\|_{\tau_{\nu}}\nonumber \\
&\le&C \|p^{n+1}-i_h p^{n+1}\|_{\tau_{\nu}^{-1}}\|\sigma_h^*(\nabla \be_h^{n+1})\|_{\tau_{\nu}} \nonumber\\
&\le& C\|p^{n+1}-i_h p^{n+1}\|_0^2+\frac{1}{4}\|\sigma_h^*(\nabla  \be_h^{n+1})\|_{\tau_{\nu}}^2\nonumber \\
&\le& C h^{2s}\|p\|_{L^\infty(H^{s})}^2 + \frac{1}{4}S_h(\be_h^{n+1},\be_h^{n+1}).
\end{eqnarray}

Collecting all the estimates we reach
\begin{eqnarray*}
\lefteqn{\|\be_h^{n+1}\|_0^2-\|\be_h^n\|_0^2+\Delta t \nu\|\nabla \be_h^{n+1}\|_0^2+\Delta t\|\sigma_h^*(\nabla \lambda_h^{n+1})\|_{\tau_p}^2
+\|\sigma_h^*(\nabla \be_h^{n+1})\|_{\tau_{\nu}}^2}\\
&\le & C \Delta t \left(1+\|\nabla \hat \bu_h^{n+1}\|_{L^\infty}+{\|\hat\bu_h^{n+1}\|_{L^\infty}^2}\right)\|\be_h^{n+1}\|_0^2
+\Delta t\|{\boldsymbol\xi}_{v_h}^{n+1}\|_0^2\nonumber\\
&&+C \Delta t h^{2s}\left((1+\nu)\|\bu\|_{L^\infty(H^{s+1})}^2+\|p\|_{L^\infty(H^s)}^2\right).\nonumber
\end{eqnarray*}
From (\ref{eq:cotahatu}) and (\ref{eq:cotanablahatu}) we deduce
\begin{equation}\label{eq:tilde_M}
1+\|\nabla \hat \bu_h^{j}\|_{L^\infty}+\|\hat\bu_h^{j}\|_{L^\infty}^2\le \tilde M_u, \quad
\tilde M_u=1+C\bigl(\|\bu\|_{L^\infty(H^{3})}+\|\bu\|_{L^\infty(H^{2})}^2\bigr).
\end{equation}
Summing up the terms, assuming that
\begin{equation}\label{eq:assume_Delta_t_tilde}
\Delta t \tilde M_u\le \frac{1}{2},
\end{equation}
and applying Lemma~\ref{gronwall} (Gronwall) leads to
\begin{eqnarray}\label{eq:after_gron"}
\lefteqn{\|\be_h^{n}\|_0^2+\Delta t\nu\sum_{j=1}^n\|\nabla \be_h^j\|_0^2+\Delta t\sum_{j=1}^n\|\sigma_h^*(\nabla \lambda_h^j)\|_{\tau_p}^2
+\Delta t\sum_{j=1}^n\|\sigma_h^*(\nabla \be_h^{j})\|_{\tau_{\nu}}^2} \\
&\le& e^{2T\tilde M_u}\left( \|\be_h^0\|_0^2+\Delta t\sum_{j=1}^n\|{\boldsymbol\xi}_{v_h}^j\|_0^2+C T h^{2s}\left((1+\nu)\|\bu\|_{L^\infty(H^{s+1})}^2+\|p\|_{L^\infty(H^s)}^2\right)\right).\nonumber
\end{eqnarray}
Now, we can argue exactly as in Section~\ref{se:velo_grad-div}
to conclude
 \begin{eqnarray}\label{eq:prinerror2}
\lefteqn{\|\be_h^{n}\|_0^2+\Delta t\nu\sum_{j=1}^n\|\nabla \be_h^j\|_0^2+\Delta t\sum_{j=1}^n\|\sigma_h^*(\nabla \lambda_h^j)\|_{\tau_p}^2
+\Delta t\sum_{j=1}^n\|\sigma_h^*(\nabla \be_h^{j})\|_{\tau_{\nu}}^2}\nonumber\\
&\le& e^{2T\tilde M_u}\left( \|\be_h^0\|_0^2+C T \tilde K_{u,p} h^{2s}+C(\Delta t)^2\int_{t_0}^{t_n}\|\partial_{tt}\bu\|_{0}^2\right),
\end{eqnarray}
with
\begin{equation}\label{eq:tilde_K}
\tilde K_{u,p}=\left((1+\|\bu\|_{L^\infty(H^{2})}^2+\nu)\|\bu\|_{L^\infty(H^{s+1})}^2+\|\partial_t\bu\|_{L^\infty(H^{s})}^2+\|p\|_{L^\infty(H^s)}^2\right).
\end{equation}
The triangle inequality finishes the proof of the velocity error estimate.

\begin{Theorem}\label{thm:velo_bound_4}
Let the solution of \eqref{eq:NSweak} be sufficiently smooth in space and time, let the
time step be sufficiently small such that \eqref{eq:assume_Delta_t_tilde} holds, and
let Assumption~A1 be satisfied. Then, the
following error bound holds for $2\le s\le l$
 \begin{eqnarray}\label{eq:err_after_gron_22}
\lefteqn{\|\bu^n-\bu_h^{n}\|_0^2+\Delta t\nu\sum_{j=1}^n\|\nabla(\bu^j - \bu_h^j)\|_0^2+\Delta t\sum_{j=1}^n\|\sigma_h^*(\nabla (p^j- p_h^j))\|_{\tau_p}^2
}\nonumber\\
&& +\Delta t\sum_{j=1}^n \|\sigma_h^*(\nabla (\bu^j - \bu_h^{j}))\|_{\tau_{\nu}}^2 \\
&\le& C e^{2T\tilde M_u}\left( \|\be_h^0\|_0^2+  T \tilde K_{u,p} h^{2s}+ (\Delta t)^2\int_{t_0}^{t_n}\|\partial_{tt}\bu\|_{0}^2~dt\right),\nonumber
\end{eqnarray}
where the constants on the right-hand side are defined in \eqref{eq:tilde_M} and \eqref{eq:tilde_K}.
\end{Theorem}

\subsection{Error bound for the pressure}
\label{sect:4pressure}

The bound for the pressure follows the lines of Section~\ref{sec:pres_bound_0} with the exception of the bound of the nonlinear term that can be handled as follows
\begin{eqnarray*}
\|B(\bu^n,\bu^n)-B(\bu_h^n,\bu_h^n)\|_{-1}\le \sup_{\|\phi\|_{1}=1}|b(\bu^n,\bu^n-\bu_h^n,\phi)|
+\sup_{\|\phi\|_{1}=1}|b(\bu^n-\bu_h^n,\bu_h^n,\phi)|.
\end{eqnarray*}
Arguing as before and recalling that $\nabla\cdot\bu =0$, we can prove
\begin{equation*}
\|B(\bu^n,\bu^n)-B(\bu_h^n,\bu_h^n)\|_{-1}
\le \left(\|\bu^n\|_{L^\infty}+\|\bu_h^n \|_{L^\infty}\right)\|\bu^n-\bu_h^n\|_0
+\sup_{\|\phi\|_{1}=1}|((\nabla \cdot \bu_h^n)\phi,\bu_h^n)|.
\end{equation*}
The last term can be decomposed as follows
\begin{equation}\label{eq:sect_42_pres_2}
((\nabla \cdot \bu_h^n)\phi,\bu_h^n)=(\nabla \cdot \bu_h^n,{\sigma_h^l}(\phi\cdot\bu_h^n))
+(\nabla \cdot \bu_h^n,(I-{\sigma_h^l})(\phi\cdot\bu_h^n)) .
\end{equation}
Since ${\sigma_h^l}(\phi\cdot\bu_h^n)) \in {\boldsymbol Y}_h^{l}$, one can use the error equation \eqref{eq:er}
for estimating the first term in \eqref{eq:sect_42_pres_2}. Applying in addition the definition
\eqref{eq:norm_tau_p} of $\|\cdot\|_{\tau_p}$, the choice \eqref{eq:tau_p} of the stabilization parameter,
the stability \eqref{eq:stasigma} of the projection, and the inverse inequality \eqref {inv} yields
\begin{eqnarray}\label{eq:pre1}
(\nabla\cdot \bu_h^n,{\sigma_h^l}(\phi\cdot\bu_h^{n}))
&\le&|s_{\rm pres}(\lambda_h^{n},{\sigma_h^l}(\phi\cdot\bu_h^{n})) |+|s_{\rm pres}(\hat p_h^{n},{\sigma_h^l}(\phi\cdot\bu_h^{n}))|
\\
&\le&Ch \left(\|\sigma^*_h(\nabla \lambda_h^{n})\|_{\tau_p}+\|\sigma^*_h(\nabla \hat p_h^{n})\|_{\tau_p}\right)
\|\sigma^*_h(\nabla {\sigma_h^l}(\phi\cdot\bu_h^n))\|_0\nonumber\\
&\le&C \left(\|\sigma^*_h(\nabla \lambda_h^{n})\|_{\tau_p}+\|\sigma^*_h(\nabla \hat p_h^{n})\|_{\tau_p}\right)
\|\phi\cdot\bu_h^n\|_0.\nonumber
\end{eqnarray}
Applying H\"older's and Sobolev's inequality, we have
$$
\|\phi\cdot\bu_h^n\|_0\le \|\phi\|_{L^{2d}}\|\bu_h^n\|_{L^{2d/(d-1)}}\le C\|\phi\|_{1} \|\bu_h^n\|_{L^{2d/(d-1)}},
$$
so that
\begin{equation*}
\sup_{\|\phi\|_{1}=1}(\nabla \cdot\bu_h^n,{\sigma_h^l}(\phi\cdot\bu_h^{n}))\le C\|\bu_h^n\|_{L^{2d/(d-1)}}\left(
\|\sigma^*_h(\nabla \lambda_h^{n})\|_{\tau_p}+C h^s\|p\|_{L^\infty(H^s)}\right).
\end{equation*}
With the decomposition
\begin{equation}\label{eq:decom_sim}
\bu_h^n - \bu^n= \be_h^n + \hat\bu_h^n-\bu^n,
\end{equation}
the inverse estimate (\ref{inv}), (\ref{eq:hatu}), and \eqref{eq:cota_sigma},
one obtains for the second term on the right-hand side of \eqref{eq:sect_42_pres_2}
\begin{equation}\label{eq:pres_est_2_02}
(\nabla \cdot \bu_h^n,(I-{\sigma_h^l})(\phi\cdot\bu_h^n))\le C h\left(h^{-1}\|\be_h^n\|_0+
h^{s}\|\bu^n\|_{s+1}\right) |\bu_h^n \phi|_{1}.
\end{equation}
The product rule and a Sobolev embedding gives
\begin{eqnarray*}
|\bu_h^n \cdot\phi|_{1}&\le& C\left( \|\bu_h^n\|_{L^\infty}|\phi|_{1}+\|\nabla \bu_h^n\|_{L^{2d/(d-1)}}\|\phi\|_{L^{2d}}\right)\\
&\le& C\left(\|\bu_h^n\|_{L^\infty}+\|\nabla \bu_h^n\|_{L^{2d/(d-1)}}\right)\|\phi\|_{1}.
\end{eqnarray*}
Now, adding and subtracting $\bu^n$, using decomposition (\ref{eq:decom_sim}) and applying the inverse inequality (\ref{inv}), (\ref{eq:prinerror2}), (\ref{eq:hatu}), and a Sobolev embedding we get
\begin{equation*}
\|\nabla \bu_h^n\|_{L^{2d/(d-1)}}\le C\left[ \frac{e^{T\tilde M_u}}{h^{3/2}}\left( \|\be_h^0\|_0^2
+ T \tilde K_{u,p} h^{2s}+(\Delta t)^2\int_{t_0}^{t_n}\|\partial_{tt}\bu\|_{0}^2\right)^{1/2}+\|\bu\|_{L^\infty(H^2)}\right].
\end{equation*}
Assuming that $s\ge 3/2$,
\begin{equation}\label{eq:sec_pres_1_ass_1}
\|\be_h^0\|_0 = \mathcal O(h^{3/2})\quad \mbox{and}\quad \Delta t\le C h^{3/2}
\end{equation}
gives
$
\|\nabla \bu_h^n\|_{L^{2d/(d-1)}}\le \tilde L_u$, where
\begin{equation}\label{eq:tilde_L_u}
\tilde L_u = C e^{T\tilde M_u} \left(\|\bu\|_{L^\infty(H^2)}^2
+ T \tilde K_{u,p}+\int_{t_0}^{t_n}\|\partial_{tt}\bu\|_{0}^2\right)^{1/2}+C\|\bu\|_{L^\infty(H^2)}.
\end{equation}
Arguing as in (\ref{eq:norma_inf_uh}), it follows that $\|\bu_h^n\|_{L^\infty}\le \tilde L_u$ whenever
$
\|\be_h^0\|_0 = \mathcal O(h^{d/2})$ and
$\Delta t  \le C h^{d/2}$,
which coincides with (\ref{eq:sec_pres_1_ass_1}) in the case $d=3$ and is weaker than (\ref{eq:sec_pres_1_ass_1}) in the case $d=2$.
Inserting the estimates  in \eqref{eq:pres_est_2_02} leads to
$$
\sup_{\|\phi\|_{1}=1}(\nabla \cdot \bu_h^n,(I-\sigma_h^l)(\phi\cdot\bu_h^n))\le \tilde L_u\left(\|\be_h^n\|_0
+h^{s+1}\|\bu^n\|_{s+1}\right).
$$
Collecting all estimates and taking into account that $\|\bu^n-\bu_h^n\|_0\le \|\be_h^n\|_0+C h^{s+1}\|\bu^n\|_{s+1}$ yields
\begin{eqnarray*}
\lefteqn{
\|B(\bu^n,\bu^n)-B(\bu_h^n,\bu_h^n)\|_{-1}}\\
&\le& \tilde L_u \big[\|\be_h^n\|_0+\|\sigma^*_h(\nabla \lambda_h^{n})\|_{\tau_p} +  h^s\left(\|p\|_{L^\infty(H^s)}+h\|\bu\|_{L^\infty(H^{s+1})}\right)\big],\nonumber
\end{eqnarray*}
and using (\ref{eq:cota_gen}) gives
\begin{eqnarray*}
\lefteqn{\sum_{j=0}^n\Delta t \|B( \bu^j,\bu^j)-B(\bu_h^j,\bu_h^j)\|_{-1}}
\nonumber\\
&\le& \tilde L_u\Bigg[T\left(\max_{1\le j\le n}\|\be_h^j\|_0+h^s\left(\|p\|_{L^\infty(H^s)}+
h\|\bu\|_{L^\infty(H^{s+1})}\right)\right)\nonumber\\
&& +T^{1/2}\left(\sum_{j=1}^n\Delta t \|\sigma^*_h(\nabla \lambda_h^{j})\|_{\tau_p}^2\right)^{1/2}\Bigg].
\end{eqnarray*}
Now, the bound for the pressure concludes as the bound of Section~\ref{sec:pres_bound_0}.

\begin{Theorem} \label{thm:pres_bound_4}
Let the assumption of Theorem~\ref{thm:velo_bound_4} and condition
\eqref{eq:sec_pres_1_ass_1}   be satisfied, then it holds
\begin{equation*}
\left\|\Delta t \sum_{j=1}^n (p^j-p_h^j)\right\|_0
\le \beta_0C(\bu,\partial_t\bu,\partial_{tt}\bu, p,T) \left(\|\bu_0 - \bu_h^0\|_0+h^s+\Delta t\right).
\end{equation*}
\end{Theorem}

\section{Local projection stabilization with control of the fluctuation of the divergence}
\label{sec:lps_grad_div}

In this section, a LPS method is briefly studied, under the same assumptions as in
Section~\ref{sec:lps_control_grad}, that uses instead of the stabilizing term
\eqref{eq:lps_gradient} a corresponding term with the divergence
\begin{equation}\label{eq:lps_grad_div}
S_h(\bu_h,\bv_h) := \sum_{K\in \mathcal{T}_h}
\tau_{\mu,K}\left(\sigma_h^* (\nabla\cdot  \bu_h), \sigma_h^*(\nabla\cdot \bv_h)\right)_K,
\end{equation}
with $\tau_{\mu,K}\sim 1$, i.e., a local projection stabilization of the grad-div term
is applied.

In Section~\ref{sec:lps_control_grad}, the stabilization with respect to the velocity
enters the error analysis in \eqref{eq:erdiv_4.1} and \eqref{eq:last_term_4.1}.
It can be readily checked that an estimate of form  \eqref{eq:erdiv_4.1} can be derived
also for \eqref{eq:lps_grad_div}. With respect to the other term, one applies
similar steps as for deriving \eqref{eq:last_term_4.1} to obtain
\begin{eqnarray*}
(\hat \lambda_h^{n+1},\nabla \cdot \be_h^{n+1})
& \le & \|p^{n+1}-i_h p^{n+1}\|_{\tau_{\mu}^{-1}}\|\sigma_h^*(\nabla \cdot \be_h^{n+1})\|_{\tau_{\mu}} \\
&\le& C\|p^{n+1}-i_h p^{n+1}\|_0^2+\frac{1}{4}\|\sigma_h^*(\nabla \cdot  \be_h^{n+1})\|_{\tau_{\mu}}^2.
\end{eqnarray*}

Altogether, the formulations of Theorems~\ref{thm:velo_bound_4} and~\ref{thm:pres_bound_4}
apply literally also to the LPS method with the local grad-div stabilization \eqref{eq:lps_grad_div}.

\begin{remark}
Let us observe that assuming $p\in H^{s+1}(\Omega)$ instead of $p\in H^s(\Omega)$ we can write
\begin{eqnarray}\label{eq:other_p}
(\hat \lambda_h^{n+1},\nabla \cdot \be_h^{n+1})&=&-(\nabla\hat \lambda_h^{n+1}, \be_h^{n+1})\nonumber\\
&\le& \|\hat \lambda_h^{n+1}\|_1\|\be_h^{n+1}\|_0
\end{eqnarray}
and then the first term is $\mathcal O(h^s)$ for $p\in H^{s+1}(\Omega)$ and the second one goes to the Gronwall lemma. This means that for equal order elements only the stabilization of the pressure gives the same rate of convergence as, for example, Galerkin plus grad-div, assuming enough regularity for the pressure.

Let us also observe that assuming $p\in H^{s+1}(\Omega)$ for the method of Section~\ref{sec:global_grad_div}, i.e., global grad-div stabilization plus LPS stabilization for the pressure, one can argue as in Section~\ref{sec:lps_control_grad} and then  apply (\ref{eq:ernonli2}) instead of (\ref{eq:ernonli}). Then, applying (\ref{eq:other_p}) instead of (\ref{eq:cotalambda}) the factor $\mu^{-1}$ disappears from  (\ref{eq:after_gron}). As a consequence,  $\mu\sim \mathcal O(h)$ is a possible option for the stabilization parameter since  with this choice
(\ref{eq:after_gron_2}) holds with $\hat K_{u,p}$ independent of $\mu^{-1}$. Let us finally point out that in view of (\ref{eq:after_gron_2})
the choice $\mu\sim \mathcal O(h)$ compared with $\mu \sim \mathcal O(1)$ gives the same rate of convergence for the $L^2$ norm of the velocity error but reduces the rate
of convergence for the divergence  by half an order.

\end{remark}
\section{A method with rate of decay $s+1/2$  of the velocity error for $\nu \le h$}\label{sec:ordersplusonehalf}

This section considers the method from Section~\ref{sec:lps_control_grad}, which adds
a stabilization term
that gives control over the fluctuation of the gradient of the velocity
and the standard LPS term for the pressure in the situation
that $\nu \le h$. It is shown that with a different choice of the stabilization
parameters and by assuming a higher regularity of the solution, both
issues compared with Section~\ref{sec:lps_control_grad}, the rate of the error decay for
the left-hand side of \eqref{eq:err_after_gron_22} can be increased to $s+1/2$.

We follow the analysis of Section~\ref{sec:lps_control_grad}.
Instead of choosing the LPS parameter for the pressure as in (\ref{eq:tau_p}),
it will be assumed that
\begin{equation}\label{eq:tau_p2}
\alpha_1h_K\le \tau_{p,K}\le \alpha_2 h_K,
\end{equation}
and instead of taking $\tau_{\nu,K}\sim 1$, it will be assumed that
\begin{equation}\label{eq:mu_2}
c_1h_K\le \tau_{\nu,K}\le c_2 h_K,
\end{equation}
with nonnegative constants $\alpha_1, \alpha_2, c_1, c_2$.
In the sequel, the assumptions for the spatial regularity of the solutions are
$p\in H^{s+1}(\Omega)$ and $\bu, \partial_t\bu \in H^{s+1}(\Omega)^d$ at almost
every time for $s\ge 2$.


The analysis starts with a different estimate of the
truncation error ${\boldsymbol \xi}_{v_h}^{n+1}$, defined in \eqref{eq:trun1}--\eqref{eq:trun12}.
In \eqref{eq:unique}, the estimate of the term coming from this error is replaced
by
$$
\|{\boldsymbol \xi}_{v_h,1}^{n+1}\|_0^2+\frac{\|\be_h^{n+1}\|_0^2}{4}+({\boldsymbol \xi}_{v_h,2}^{n+1},\be_h^{n+1}).
$$
The term
$({\boldsymbol \xi}_{v_h,2}^{n+1},\be_h^{n+1})$ can be decomposed in the form
\begin{eqnarray}\label{eq:loscuatro}
\lefteqn{
|b(\bu^{n+1},\bu^{n+1},\be_h^{n+1})-b(\hat\bu_h^{n+1}, \hat \bu_h^{n+1},\be_h^{n+1})|}
\nonumber\\
&\le&|((\hat\bu_h^{n+1}\cdot \nabla) (\hat\bu_h^{n+1}-\bu^{n+1}),\be_h^{n+1})|
+\frac{1}{2}|((\nabla \cdot\hat\bu_h^{n+1})(\hat\bu_h^{n+1}-\bu^{n+1}),\be_h^{n+1})|
\nonumber\\
&&+|(((\hat \bu_h^{n+1}-\bu^{n+1})\cdot \nabla) \bu^{n+1},\be_h^{n+1})+\frac{1}{2}|(\nabla \cdot (\hat \bu_h^{n+1}-\bu^{n+1})\bu^{n+1},\be_h^{n+1})|.
\end{eqnarray}
Since $\|\nabla \bu^{n+1}\|_{L^\infty}$ is bounded by the regularity assumption and
$\|\nabla \cdot \hat \bu_h^{n+1}\|_{L^\infty}$ is bounded in
\eqref{eq:cotanablahatu}, the second and third terms in (\ref{eq:loscuatro}) can be bounded by
$$C\|\bu\|_{L^\infty(H^3)}\|\hat\bu_h^{n+1}-\bu^{n+1}\|_0\|\be_h^{n+1}\|_0.$$ Thus, we only need to bound the first and the last term in (\ref{eq:loscuatro}). Using integration by parts gives the decomposition
\begin{eqnarray*}
(\hat\bu_h^{n+1}\cdot \nabla (\hat\bu_h^{n+1}-\bu^{n+1}),\be_h^{n+1})
&=&-((\nabla \cdot \hat \bu_h^{n+1})(\hat \bu_h^{n+1}-\bu^{n+1}),\be_h^{n+1})\\
&&-(\hat \bu_h^{n+1}\cdot \nabla \be_h^{n+1},\hat \bu_h^{n+1}-\bu^{n+1}).
\end{eqnarray*}
Again, the first term can be bounded by $C\|\bu\|_{L^\infty(H^3)}\|\hat\bu_h^{n+1}-\bu^{n+1}\|_0\|\be_h^{n+1}\|_0$, so we only need to bound the second one. Using that the range
of $\sigma_h^{l-1}$ is ${\boldsymbol Y}_h^{l-1}$ and the definition \eqref{eq:orto} of
$\hat \bu_h^{n+1}$ yields
\begin{eqnarray}\label{eq:truncaul}
\lefteqn{(\hat \bu_h^{n+1}\cdot \nabla \be_h^{n+1},\hat \bu_h^{n+1}-\bu^{n+1})}\nonumber\\
&=&(\sigma_h^*(\hat \bu_h^{n+1}\cdot \nabla \be_h^{n+1}),\hat\bu_h^{n+1}-\bu^{n+1})\\
&=&
(\sigma_h^*(\hat \bu_h^{n+1}\cdot \sigma_h^{l-1}\nabla \be_h^{n+1}),\hat \bu_h^{n+1}-\bu^{n+1})+(\sigma_h^*(\hat \bu_h^{n+1}\cdot \sigma_h^*\nabla \be_h^{n+1}),\hat \bu_h^{n+1}-\bu^{n+1}).\nonumber
\end{eqnarray}
We apply Lemma~\ref{cor_burman} to the first term to obtain
\begin{eqnarray*}
\lefteqn{|(\sigma_h^*(\hat \bu_h^{n+1}\cdot \sigma_h^{l-1}\nabla \be_h^{n+1}),\hat \bu_h^{n+1}-\bu^{n+1})|}\nonumber\\
&\le& C h \|\hat \bu_h^{n+1}\|_{W^{1,\infty}}\|\sigma_h^{l-1}\nabla \be_h^{n+1}\|_0\|\hat \bu_h^{n+1}-\bu^{n+1}\|_0\nonumber\\
&\le& C\|\hat \bu_h^{n+1}\|_{W^{1,\infty}} \|\be_h^{n+1}\|_0\|\hat \bu_h^{n+1}-\bu^{n+1}\|_0,
\end{eqnarray*}
where in the last inequality we have applied the~$L^2$ stability of~$\sigma_h^{l-1}$
\eqref{eq:stasigma} and the inverse
inequality (\ref{inv}).
For the second term of \eqref{eq:truncaul}, we get with \eqref{eq:stasigma}
\begin{eqnarray*}
\lefteqn{|(\sigma_h^*(\hat \bu_h^{n+1}\cdot \sigma_h^*\nabla \be_h^{n+1}),\hat \bu_h^{n+1}-\bu^{n+1})|}
\nonumber\\
&\le & C
\sum_{K\in {\mathcal T}_h}\|\hat \bu_h^{n+1}\cdot \sigma_h^*\nabla \be_h^{n+1}\|_{L^2(\omega_K)}\|\hat \bu_h^{n+1}-\bu^{n+1}\|_{L^2(K)}
\nonumber\\
&\le & C
\sum_{K\in {\mathcal T}_h}\|\hat \bu_h^{n+1}\cdot \sigma_h^*\nabla \be_h^{n+1}\|_{L^2(K)}\|\hat \bu_h^{n+1}-\bu^{n+1}\|_{L^2(K)}
\nonumber\\
&\le& C
\sum_{K\in {\mathcal T}_h}\|\hat \bu_h^{n+1}\|_{L^{\infty}(K)}\| \sigma_h^*\nabla \be_h^{n+1}\|_{L^2(K)}\|\hat \bu_h^{n+1}-\bu^{n+1}\|_{L^2(K)}
\nonumber\\
&\le& C\|\hat \bu_h^{n+1}\|_{L^\infty}^2\sum_{K\in {\mathcal T}_h}\tau_{\nu,K}^{-1}\|\hat \bu_h^{n+1}-\bu^{n+1}\|_{L^2(K)}^2
+\frac{1}{8}\sum_{K\in {\mathcal T}_h}\tau_{\nu,K}\|\sigma_h^*\nabla \be_h^{n+1}\|_{L^2(K)}^2.
\end{eqnarray*}
This bound concludes the estimate of the first term on the right-hand side of (\ref{eq:loscuatro}).
To bound the last term on the right-hand side of (\ref{eq:loscuatro}), integration by parts and
\eqref{eq:orto} are applied
\begin{eqnarray*}
\lefteqn{
|\nabla \cdot (\hat \bu^{n+1}_h-\bu^{n+1})\bu^{n+1},\be_h^{n+1})|} \nonumber\\
&=&|-(\hat \bu^{n+1}_h-\bu^{n+1},\sigma_h^*\nabla(\bu^{n+1}\cdot\be_h^{n+1}))|\nonumber\\
&\le&|(\hat \bu^{n+1}_h-\bu^{n+1},\sigma_h^*(\nabla \bu^{n+1}\be_h^{n+1}))|
+|(\hat \bu^{n+1}_h-\bu^{n+1},\sigma_h^*(\nabla \be_h^{n+1}\bu^{n+1}))|\nonumber\\
&\le&\|\hat \bu^{n+1}_h-\bu^{n+1}\|_0\|\nabla \bu^{n+1}\|_{L^\infty}\|\be_h^{n+1}\|_0
+|(\hat \bu^{n+1}_h-\bu^{n+1},\sigma_h^*(\nabla \be_h^{n+1}\bu^{n+1}))|.
\end{eqnarray*}
The last term can be bounded arguing exactly as in (\ref{eq:truncaul}). Thus, collecting
all estimates and using \eqref{eq:cotahatu} to bound $\|\hat\bu_h^{n+1}\|_{L^\infty} \le C\|\bu\|_{L^\infty(H^2)} $  yields
\begin{eqnarray}
\lefteqn{
|b(\bu^{n+1},\bu^{n+1},\be_h^{n+1})-b(\hat\bu_h^{n+1},\hat \bu_h^{n+1},\be_h^{n+1})|}\nonumber\\
&\le &  C\|\bu\|_{L^\infty(H^3)}\|\bu^{n+1}-\hat \bu_h^{n+1}\|_0\|\be_h^{n+1}\|_0
+ C\|\bu\|_{L^\infty(H^2)}^2\sum_{K\in {\mathcal T}_h}\tau_{\nu,K}^{-1}\|\hat \bu_h^{n+1}-\bu^{n+1}\|_{L^2(K)}^2\nonumber\\
&& +\frac{1}{4}\sum_{K\in {\mathcal T}_h}\tau_{\nu,K}\|\sigma_h^*\nabla \be_h^{n+1}\|_{L^2(K)}^2\nonumber\\
&\le &
C\|\bu\|_{L^\infty(H^3)}^2\|\bu^{n+1}-\hat \bu_h^{n+1}\|_0^2+\frac{1}{4}\|\be_h^{n+1}\|_0^2
+ C\|\bu\|_{L^\infty(H^2)}^2\sum_{K\in {\mathcal T}_h}\tau_{\nu,K}^{-1}\|\hat \bu_h^{n+1}-\bu^{n+1}\|_{L^2(K)}^2 \nonumber\\
&&
+\frac{1}{4}\sum_{K\in {\mathcal T}_h}\tau_{\nu,K}\|\sigma_h^*\nabla \be_h^{n+1}\|_{L^2(K)}^2\nonumber
\\
&\le & C\|\bu\|_{L^\infty(H^3)}^2\left(\max_{K\in {\mathcal T}_h}\tau_{\nu,K}^{-1}\right)\|\bu^{n+1}-\hat\bu_h^{n+1}\|_{0}^2
+\frac{1}{4}\|\be_h^{n+1}\|_0^2+\frac{1}{4}S_h(\be_h^{n+1},\be_h^{n+1}),
\end{eqnarray}
where we have bounded
$\min_{K\in {\mathcal T}_h}\{\tau_{\nu,K}\}\|\bu\|_{L^\infty(H^3)}^2 +
\|\bu\|_{L^\infty(H^2)}^2\le C \|\bu\|_{L^\infty(H^3)}^2 +
\|\bu\|_{L^\infty(H^2)}^2\le C \|\bu\|_{L^\infty(H^3)}^2$.

Thus, in the present case, instead of~(\ref{eq:unique}), we have
\begin{eqnarray}\label{eq:unique2}
\lefteqn{\frac{\|\be_h^{n+1}\|_0^2}{2\Delta t}-\frac{\|\be_h^{n}\|_0^2}{2\Delta t}+\frac{\|\be_h^{n+1}-\be_h^n\|_0^2}{2\Delta t}
+\frac{\nu}{2}\|\nabla \be_h^{n+1}\|_0^2
+\|\sigma_h^*(\nabla \lambda_h^{n+1})\|_{\tau_p}^2}
\nonumber\\
&&+\frac{3}{4}S_h( \be_h^{n+1}, \be_h^{n+1})\nonumber\\
&\le& \left|b(\bu_h^{n+1},\bu_h^{n+1},\be_h^{n+1})-b(\hat \bu_h^{n+1},\hat \bu_h^{n+1},\be_h^{n+1})\right|+\|{\boldsymbol \xi}_{v_{h,1}}^{n+1}\|_0^2+\frac{\|\be_h^{n+1}\|_0^2}{2}\\
&&+\left|({\boldsymbol \xi}_{q_h}^{n+1},\lambda_h^{n+1})\right|
+\frac{\nu}{2}\|\nabla\hat \be_h^{n+1}\|_0
+\left|s_{\rm  pres}(\hat p_h^{n+1},\lambda_h^{n+1})\right|\nonumber\\
&&+C\|\bu\|_{L^\infty(H^3)}^2\left(\max_{K\in {\mathcal T}_h}\tau_{\nu,K}^{-1}\right)\|\hat\be_h^{n+1}\|_0^2
+\left|S_h(\hat\bu_h^{n+1},\be_h^{n+1})\right|
+\left|(\hat \lambda_h^{n+1},\nabla \cdot \be_h)\right|.\nonumber
\end{eqnarray}
Next, we argue as in Section~\ref{sec:lps_control_grad} and apply (\ref{eq:esta}),
(\ref{eq:likelu1}), and (\ref{eq:52}) as starting point for estimating the first term on the right-hand side
of \eqref{eq:unique2}.
To bound the first term on the right-hand side of  (\ref{eq:52}), a similar approach as in \eqref{eq:handle00} is applied, taking into account
the different stabilization parameter and regularity of the solution,
\begin{eqnarray}\label{eq:handle}
\lefteqn{s_{\rm pres}(p_h^{n+1},\sigma_h^l(\hat \bu_h^{n+1}\cdot\be_h^{n+1}))}\\
&\le& C h^{2s+1}\|p\|_{L^\infty(H^{s+1})}^2+\frac{1}{8}\|\sigma^*_h(\nabla \lambda_h^{n+1})\|_{\tau_p}^2+4\left( \max_{K\in\mathcal T_h} \tau_{p,K}\right)\|\sigma^*_h(\nabla \sigma_h^l(\hat \bu_h^{n+1}\cdot\be_h^{n+1}))\|_0^2.\nonumber
\end{eqnarray}
Now, the bound of the last term of \eqref{eq:handle} becomes different as in Section~\ref{sec:lps_control_grad}
since the application of the inverse inequality gives rise to a term with factor $h^{-1}$, compare \eqref{eq:handle00}.
The triangle inequality gives
\begin{eqnarray}\label{eq:medio1}
\|\sigma^*_h(\nabla \sigma_h^l(\hat \bu_h^{n+1}\cdot\be_h^{n+1}))\|_0^2&\le& 2\|\sigma^*_h(\nabla (\hat \bu_h^{n+1}\cdot\be_h^{n+1}))\|_0^2
\nonumber\\
&&+2\|\sigma^*_h(\nabla (I-\sigma_h^l)(\hat \bu_h^{n+1}\cdot\be_h^{n+1}))\|_0^2.
\end{eqnarray}
For the second term on the right-hand side of \eqref{eq:medio1}, we apply the $L^2$ stability \eqref{eq:stasigma} of $\sigma_h^*$  and (\ref{eq:berto}) to get
\begin{eqnarray}
\|\sigma^*_h(\nabla (I-\sigma_h^l)(\hat \bu_h^{n+1}\cdot\be_h^{n+1}))\|_0^2
&\le& C\|\nabla (I-\sigma_h^l)(\hat \bu_h^{n+1}\cdot\be_h^{n+1})\|_0^2\nonumber\\
&\le& C \|\hat \bu_h^{n+1} \|_{W^{1,\infty}}^2\|\be_h^{n+1}\|_0^2.\label{eq:medio12}
\end{eqnarray}
Utilizing the product rule, the triangle inequality, and \eqref{eq:stasigma} gives
for the first term on the right-hand side of (\ref{eq:medio1})
\begin{equation}\label{eq:half3}
\|\sigma^*_h(\nabla (\hat \bu_h^{n+1}\cdot \be_h^{n+1}))\|_0\le  C \|\nabla \hat \bu_h^{n+1}\|_{L^\infty}\|\be_h^{n+1}\|_0+
\|\sigma_h^*(\nabla \be_h^{n+1}\hat \bu_h^{n+1})\|_0.
\end{equation}
For the second term on the right-hand side of \eqref{eq:half3}, we use the decomposition
$\nabla \be_h^{n+1} = \sigma_h^{l-1}\nabla\be_h^{n+1}  + \sigma^* \nabla\be_h^{n+1} $,
Lemma~\ref{cor_burman}, \eqref{eq:stasigma}, and the inverse estimate \eqref{inv}
to obtain
\begin{eqnarray}\label{eq:half4}
\|\sigma_h^*(\nabla \be_h^{n+1}\hat \bu_h^{n+1})\|_0
&\le& Ch \|\hat \bu_h^{n+1}\|_{W^{1,\infty}}\|\sigma_h^{l-1} \nabla\be_h^{n+1}\|_0+
\|\sigma_h^*((\sigma_h^*\nabla \be_h^{n+1})\hat \bu_h^{n+1})\|_0\nonumber\\
&\le& C\|\hat \bu_h^{n+1}\|_{W^{1,\infty}}\|\be_h^{n+1}\|_0+
C\|(\sigma_h^*\nabla \be_h^{n+1})\hat \bu_h^{n+1}\|_0.
\end{eqnarray}
For the second term on the right-hand-side of (\ref{eq:half4}) we get
\begin{eqnarray}\label{eq:half5}
\|(\sigma_h^*\nabla \be_h^{n+1})\hat\bu_h^{n+1}\|_0^2&=&\sum_{K\in\mathcal T_h}\|(\sigma_h^*\nabla \be_h^{n+1})\hat\bu_h^{n+1}\|_{L^2(K)}^2
\nonumber\\
&\le& \sum_{K\in\mathcal T_h}\|\hat\bu_h^{n+1}\|_{L^\infty(K)}^2\|\sigma_h^*\nabla \be_h^{n+1}\|_{L^2(K)}^2\nonumber\\
&=&\sum_{K\in\mathcal T_h}\tau_{\nu,K}^{-1}\|\hat\bu_h^{n+1}\|_{L^\infty(K)}^2\tau_{\nu,K}\|\sigma_h^*\nabla \be_h^{n+1}\|_{L^2(K)}^2
\nonumber\\
&\le&\left(\max_{K\in\mathcal T_h} \tau_{\nu,K}^{-1}\right)\|\hat\bu_h^{n+1}\|_{L^\infty}^2\|\sigma_h^*(\nabla\be_h^{n+1})\|_{\tau_\nu}^2.
\end{eqnarray}
Altogether, we conclude from (\ref{eq:half3}), (\ref{eq:half4}), and (\ref{eq:half5}) that
\begin{eqnarray}\label{eq:medio5}
\|\sigma^*_h(\nabla (\hat \bu_h^{n+1}\cdot\be_h^{n+1}))\|_0^2&\le&C\|\hat\bu_h^{n+1}\|_{W^{1,\infty}}^2\|\be_h^{n+1}\|_0^2\\
&&+C\|\hat\bu_h^{n+1}\|_{L^\infty}^2\left(\max_{K\in\mathcal T_h} \tau_{\nu,K}^{-1}\right)\|\sigma_h^*(\nabla \be_h^{n+1})\|_{\tau_\nu}^2.\nonumber
\end{eqnarray}
Taking into account (\ref{eq:medio1}), (\ref{eq:medio12}), and (\ref{eq:medio5}), we finally obtain for the last term on the right-hand side of
(\ref{eq:handle})
\begin{eqnarray}\label{eq:lastspres0}
\lefteqn{4\left( \max_{K\in\mathcal T_h} \tau_{p,K}\right)\|\sigma^*_h(\nabla \sigma_h^l(\hat \bu_h^{n+1}\cdot \be_h^{n+1}))\|_0^2}\nonumber\\
&\le&Ch\|\hat\bu_h^{n+1}\|_{W^{1,\infty}}^2\|\be_h^{n+1}\|_0^2\nonumber\\
&&
+C \|\hat\bu_h^{n+1}\|_{L^\infty}^2\left( \max_{K\in\mathcal T_h} \tau_{p,K}\right)\left(\max_{K\in\mathcal T_h} \tau_{\nu,K}^{-1}\right)\|\sigma_h^*(\nabla\be_h^{n+1})\|_{\tau_\nu}^2.
\end{eqnarray}
Thus, assuming
\begin{equation}\label{eq:cond_dif}
C \|\hat\bu_h^{n+1}\|_{L^\infty}^2\left( \max_{K\in\mathcal T_h} \tau_{p,K}\right)\left(\max_{K\in\mathcal T_h} \tau_{\nu,K}^{-1}\right)\le \frac {1}{16},
\end{equation}
with $C$ being the constant of the last term of  (\ref{eq:lastspres0}),  estimate (\ref{eq:lastspres0}) gives
\begin{eqnarray}\label{eq:lastspres}
\lefteqn{
4\left( \max_{K\in\mathcal T_h} \tau_{p,K}\right)\|\sigma^*_h(\nabla \sigma_h^l(\hat \bu_h^{n+1}\cdot \be_h^{n+1}))\|_0^2}\nonumber\\
&\le&Ch\|\hat\bu_h^{n+1}\|_{W^{1,\infty}}^2\|\be_h^{n+1}\|_0^2
+\frac{1}{16}S_h(\be_h^{n+1},\be_h^{n+1}).
\end{eqnarray}
From~(\ref{eq:handle}) and (\ref{eq:lastspres}) we get now
\begin{eqnarray}\label{eq:handle2}
s_{\rm pres}(p_h^{n+1},\sigma_h^l(\hat \bu_h^{n+1}\cdot\be_h^{n+1}))
&\le& C h^{2s+1}\|p\|_{L^\infty(H^{s+1})}^2+Ch
\left\|\hat\bu_h^{n+1}\right\|_{W^{1,\infty}}^2
\|\be_h^{n+1}\|_0^2,\nonumber\\
&&+\frac{1}{8}\|\sigma^*_h(\nabla \lambda_h^{n+1})\|_{\tau_p}^2+\frac{1}{16}S_h(\be_h^{n+1},\be_h^{n+1}).
\end{eqnarray}
Observe that~(\ref{eq:handle2}) is the counterpart of (\ref{eq:handle00}).

To bound the second term on the right-hand side of (\ref{eq:52}), applying integration
by parts, \eqref{eq:orto}, the Cauchy--Schwarz inequality, and Young's inequality yields
\begin{eqnarray}\label{eq:97a}
\left((\nabla \cdot \hat\be_h^{n+1}),\sigma_h^l(\hat \bu_h^{n+1}\cdot\be_h^{n+1})\right)&=&-\left(\hat \be_h^{n+1},\sigma^*_h(\nabla\sigma_h^l
(\hat \bu_h^{n+1}\cdot\be_h^{n+1}))\right)\nonumber\\
&\le&\frac{\|\hat\be_h^{n+1}\|_0^2}{4\varepsilon h}+\varepsilon h\|\sigma^*_h(\nabla\sigma_h^l(\hat \bu_h^{n+1}\cdot\be_h^{n+1}))\|_0^2\\
&\le&C\varepsilon^{-1}h^{2s+1}\|\bu\|_{L^\infty(H^{s+1})}+\varepsilon h\|\sigma^*_h(\nabla\sigma_h^l(\hat \bu_h^{n+1}\cdot\be_h^{n+1}))\|_0^2\nonumber
\end{eqnarray}
with some $\varepsilon > 0$. Now, the second term on the right-hand side can be estimated the
same way as the second term of \eqref{eq:handle}. The parameter $\varepsilon$ can be chosen
sufficiently small so that
\begin{equation}
\label{eq:epsilon}
C \varepsilon h \|\hat\bu_h^{n+1}\|_{L^\infty}^2\left(\max_{K\in\mathcal T_h} \tau_{\nu,K}^{-1}\right)\le \frac {1}{16},
\end{equation}
and hence, the second term of \eqref{eq:97a} can be bounded by  (\ref{eq:lastspres}).

Collecting terms and assuming that condition (\ref{eq:cond_dif}) holds, instead of (\ref{eq:ernonli2}), we reach
\begin{eqnarray*}
\lefteqn{
\left|b(\bu_h^{n+1},\bu_h^{n+1},\be_h^{n+1})-b(\hat \bu_h^{n+1},\hat \bu_h^{n+1},\be_h^{n+1})\right|}\nonumber\\
&\le&
C\left(\|\nabla \hat \bu_h^{n+1}\|_{L^\infty}+h\|\hat \bu_h^{n+1}\|_{W^{1,\infty}}^2 \right)\|\be_h^{n+1}\|_0^2
 +\frac{1}{8}\|\sigma^*_h(\nabla \lambda_h^{n+1})\|_{\tau_p}^2\nonumber\\
&&+\frac{1}{8}S_h(\be_h^{n+1},\be_h^{n+1})
+C h^{2s+1}\left(\|p\|_{L^\infty(H^{s+1})}^2+\varepsilon^{-1}\|\bu\|_{L^\infty(H^{s+1})}^2\right).
\end{eqnarray*}

Now, we argue as in Section~\ref{sec:lps_control_grad}, taking into account that $p\in H^{s+1}(\Omega)$ and applying (\ref{eq:tau_p2}) and (\ref{eq:mu_2}). The estimate of the fourth term on the right-hand side
of \eqref{eq:unique2} uses the approach of \eqref{eq:erdiv2} and the choice of the stabilization parameter
\eqref{eq:tau_p2}. The seventh term is bounded by \eqref{eq:hatu} and the stabilization parameter
\eqref{eq:mu_2}.
To get a higher order of
the fifth term of \eqref{eq:unique}, we have to assume that
\begin{equation}\label{eq:ass_sect6_1}
\nu \le h.
\end{equation}
Collecting all estimates gives, instead of \eqref{eq:after_gron"},
 \begin{eqnarray*}
&&\|\be_h^{n}\|_0^2+\Delta t\nu\sum_{j=1}^n\|\nabla \be_h^j\|_0^2+\Delta t\sum_{j=1}^n\|\sigma_h^*(\nabla \lambda_h^j)\|_{\tau_p}^2
+\frac{\Delta t}{4}\sum_{j=1}^n\|\sigma_h^*(\nabla \be_h^{j})\|_{\tau_{\nu}}^2\qquad\qquad\\
&&\quad  \le  e^{2T M_u}\left( \|\be_h^0\|_0^2+2\Delta t\sum_{j=1}^n\|{\boldsymbol\xi}_{v_h,1}^j\|_0^2+C T h^{2s+1}\left(\|\bu\|_{L^\infty(H^{s+1})}^2+\|p\|_{L^\infty(H^{s+1})}^2\right)\right),\nonumber
\end{eqnarray*}
where
\begin{equation}\label{eq:Mu}
1+C\left(\|\nabla \hat \bu_h^{n+1}\|_{L^\infty}+h\|\hat \bu_h^{n+1}\|_{W^{1,\infty}}^2 \right)\le M_u=1+C\|\bu\|_{L^\infty(H^3)}\left(1+\|\bu\|_{L^\infty(H^3)}\right).
\end{equation}
Note that  we apply (\ref{eq:trunca1_1}) and (\ref{eq:trunca1_2}) under the assumption $\partial_t \bu \in H^{s+1}(\Omega)^d$ to bound $\|{\boldsymbol\xi}_{v_h,1}^j\|_0^2$. Then, instead of (\ref{eq:prinerror2}), we obtain
\begin{eqnarray*}
\lefteqn{\|\be_h^{n}\|_0^2+\Delta t\nu\sum_{j=1}^n\nu\|\nabla \be_h^j\|_0^2+\Delta t\sum_{j=1}^n\|\sigma_h^*(\nabla \lambda_h^j)\|_{\tau_p}^2
+\Delta t\sum_{j=1}^n\|\sigma_h^*(\nabla \be_h^{j})\|_{\tau_{\nu}}^2}\nonumber\\
&\le& e^{2T M_u}\left( \|\be_h^0\|_0^2+C T K_{u,p} h^{2s+1}+C(\Delta t)^2\int_{t_0}^{t_n}\|\partial_{tt}\bu\|_{0}^2\right),
\end{eqnarray*}
with
\begin{equation}\label{eq:tilde_Kb}
K_{u,p}=\left(\left(1+\varepsilon^{-1}+\|\bu\|_{L^\infty(H^3)}^2\right)\|\bu\|_{L^\infty(H^{s+1})}^2+\|\partial_t\bu\|_{L^\infty(H^{s+1})}^2+\|p\|_{L^\infty(H^{s+1})}^2\right),
\end{equation}
$\varepsilon$ being the value in~(\ref{eq:epsilon}).
The triangle inequality finishes the proof of the velocity error estimate.

\begin{Theorem}\label{thm:velo_bound_4b}
Let the assumptions of  Theorem~\ref{thm:velo_bound_4} be satisfied, let in particular
$\bu, \partial_t\bu \in L^\infty(0,T;H^{s+1}(\Omega)^d)$ and
$p\in L^\infty(0,T;H^{s+1}(\Omega))$. Let the stabilization parameters be chosen such that
 \eqref{eq:cond_dif} is satisfied and let condition
\eqref{eq:ass_sect6_1} hold.
Then, the
following error bound is valid
 \begin{eqnarray}\label{eq:err_after_gron_22b}
\lefteqn{\|\bu^n-\bu_h^{n}\|_0^2+\Delta t\nu\sum_{j=1}^n\|\nabla(\bu^j - \bu_h^j)\|_0^2+\Delta t\sum_{j=1}^n\|\sigma_h^*(\nabla (p^j- p_h^j))\|_{\tau_p}^2
}\nonumber\\
&& +\Delta t\sum_{j=1}^n \|\sigma_h^*(\nabla (\bu^j - \bu_h^{j}))\|_{\tau_{\nu}}^2 \\
&\le& C e^{2T M_u}\left( \|\be_h^0\|_0^2+  T  K_{u,p} h^{2s+1}+ (\Delta t)^2\int_{t_0}^{t_n}\|\partial_{tt}\bu\|_{0}^2~dt\right),\nonumber
\end{eqnarray}
where the constants on the right-hand side are defined in \eqref{eq:Mu} and \eqref{eq:tilde_Kb}.
\end{Theorem}


\begin{remark}
The bound for the pressure follows the steps of Section~\ref{sect:4pressure} with the only difference that due to the change in the size of the pressure stabilization parameter instead
of (\ref{eq:pre1})  we get
\begin{eqnarray*}
(\nabla \cdot(\bu_h^n-\bu^n),\sigma_h^l(\phi\cdot\bu_h^{n}))
\le Ch^{-1/2}\left(\|\sigma^*(\nabla \lambda_h^{n})\|_{\tau_p}+\|\sigma^*_h(\nabla \hat p_h^{n})\|_{\tau_p}\right)
\|\phi\cdot\bu_h^n\|_{0},
\end{eqnarray*}
and
\begin{eqnarray*}
&&\sup_{\|\phi\|_{1}=1}(\nabla \cdot(\bu_h^n-\bu^n),\sigma_h^l(\phi\cdot\bu_h^{n}))\nonumber\\
&&\quad \quad\le C\|\bu_h^n\|_{L^{2d/(d-1)}}\left(h^{-1/2}
\|\sigma^*_h(\nabla \lambda_h^{n})\|_{\tau_p}+C h^{s+1/2}\|p\|_{L^\infty(H^{s+1})}\right).
\end{eqnarray*}
The factor $h^{-1/2}$ remains during the analysis in front of $\|\sigma^*_h(\nabla \lambda_h^{n})\|_{\tau_p}$ such that a higher rate of error decay for the pressure error cannot be proved with
this approach.

The last term in the second line of \eqref{eq:pre1} has the same principal form as the
last term of \eqref{eq:handle}. In contrast to the analysis for the velocity, we did not find a way
to replace the application of the inverse estimate by a more sophisticated approach that
leads to an improvement of the rate of error decay for the pressure.
\end{remark}

\section{Numerical studies}\label{sec:numres}

Numerical studies will be presented for the sake of supporting the analytical results.
Simulations were performed at a problem defined in $\Omega = (0,1)^2$ and the time interval
$(0,5]$ with the prescribed solution
\begin{eqnarray*}
\bu &=& \cos(t)\begin{pmatrix} \sin(\pi x-0.7)\sin(\pi y + 0.2) \\
\cos(\pi x-0.7)\cos(\pi y + 0.2)
\end{pmatrix}, \\
p &=& \cos(t)(\sin(x)\cos(y)+(\cos(1)-1)\sin(1)).
\end{eqnarray*}

The version of the Scott--Zhang operator proposed in \cite{Bad12} was used for computing the
local projection. The numerical studies were performed with the code {\sc MooNMD} \cite{JM04}.

\begin{figure}[t!]
\centerline{\includegraphics[width=0.25\textwidth]{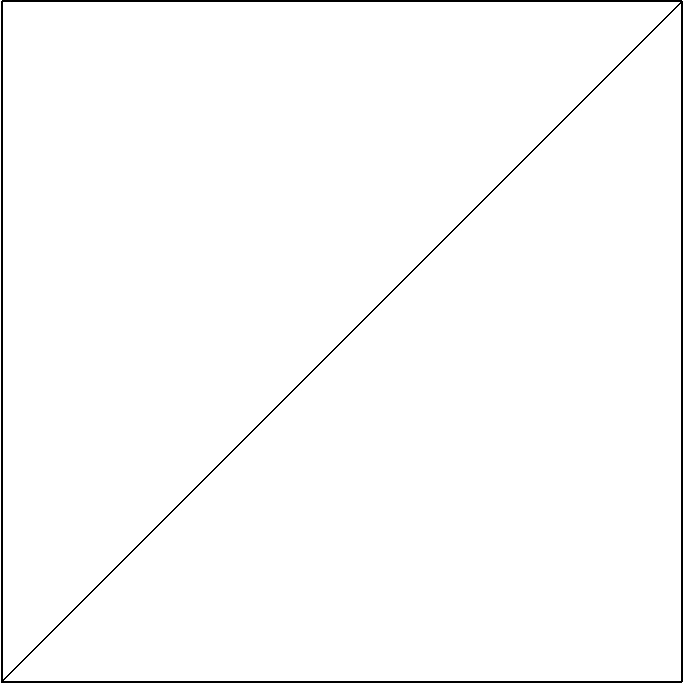}\hspace*{1em}
\includegraphics[width=0.25\textwidth]{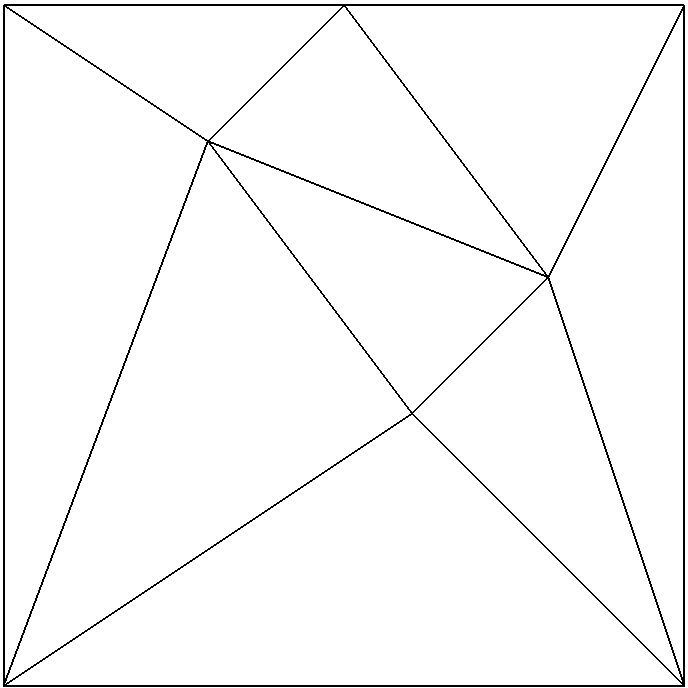}
}
\caption{Grid~1 and~2, level~0.}\label{fig:grids}
\end{figure}

We will concentrate on the convergence with respect to the mesh width. As
temporal discretization, the second order Crank--Nicolson scheme with the small
time step $\Delta t = 0.001$ was applied. Hence, the temporal error possesses a
negligible impact on the first refinements of the coarsest grids presented in
Figure~\ref{fig:grids}. The nonlinear problems in each discrete time were solved
until the Euclidean norm of the residual vector was less than $10^{-13}$.

\subsection{LPS with global grad-div stabilization}\label{sec:numres_method3}

Here, method \eqref{eq:gal} analyzed in Section~\ref{sec:global_grad_div}, with the
Crank--Nicolson scheme instead of the implicit Euler method, will be
studied.

The asymptotic choice of the LPS stabilization parameter is given in \eqref{eq:tau_p}.
From numerical studies, we could see that $\tau_{p,K}  = h_K^2$ is an appropriate
selection with respect to the accuracy of the computational results.  From the statements
of Theorem~\ref{thm:velo_bound} and~\ref{thm:pres_bound},
it follows that the grad-div stabilization parameter should be a constant.
Numerical tests showed that $\mu = 0.1$ is a good choice. In addition,
since in the considered example the pressure solution is smooth, it would be possible to
obtain in the last term of \eqref{eq:cotalambda}
$$
\frac{C}{\mu} h^{2(s+1)}\|p\|^2_{L^\infty(H^{s+1})},
$$
such that also the choice $\mu\sim h$ is possible without reducing the order of convergence.
Thus, also results for $\mu = 0.1h_K$ will be presented. Note that $\mu\sim h$ is the choice that is
proposed for the equal-order SUPG/PSPG/grad-div stabilized finite element method of the
Oseen equations, compare \cite[Rem.~5.42]{Joh16}.

Besides a number of standard errors, an error is monitored
that is
an approximation of the left-hand side of \eqref{eq:after_gron_2}. The approximation consists
in considering instead of the pressure term, the term
\begin{equation}\label{eq:press_err_sim}
\Delta t \sum_{j=1}^n \tau_p \|\nabla(p^j- p_h^j) \|_0^2,
\end{equation}
with $\tau_p = h^2$ and $h= h_0 2^{-l}$, $l$ being the index of the level
with $h_0 = \sqrt2$ for Grid~1 and $h_0= 1$ for Grid~2.  Using \eqref{eq:stasigma}, the pressure term
on the left-hand side of \eqref{eq:after_gron_2} can be estimated from above with \eqref{eq:press_err_sim}
times a constant.

\begin{figure}[t!]
\centerline{\includegraphics[width=0.4\textwidth]
{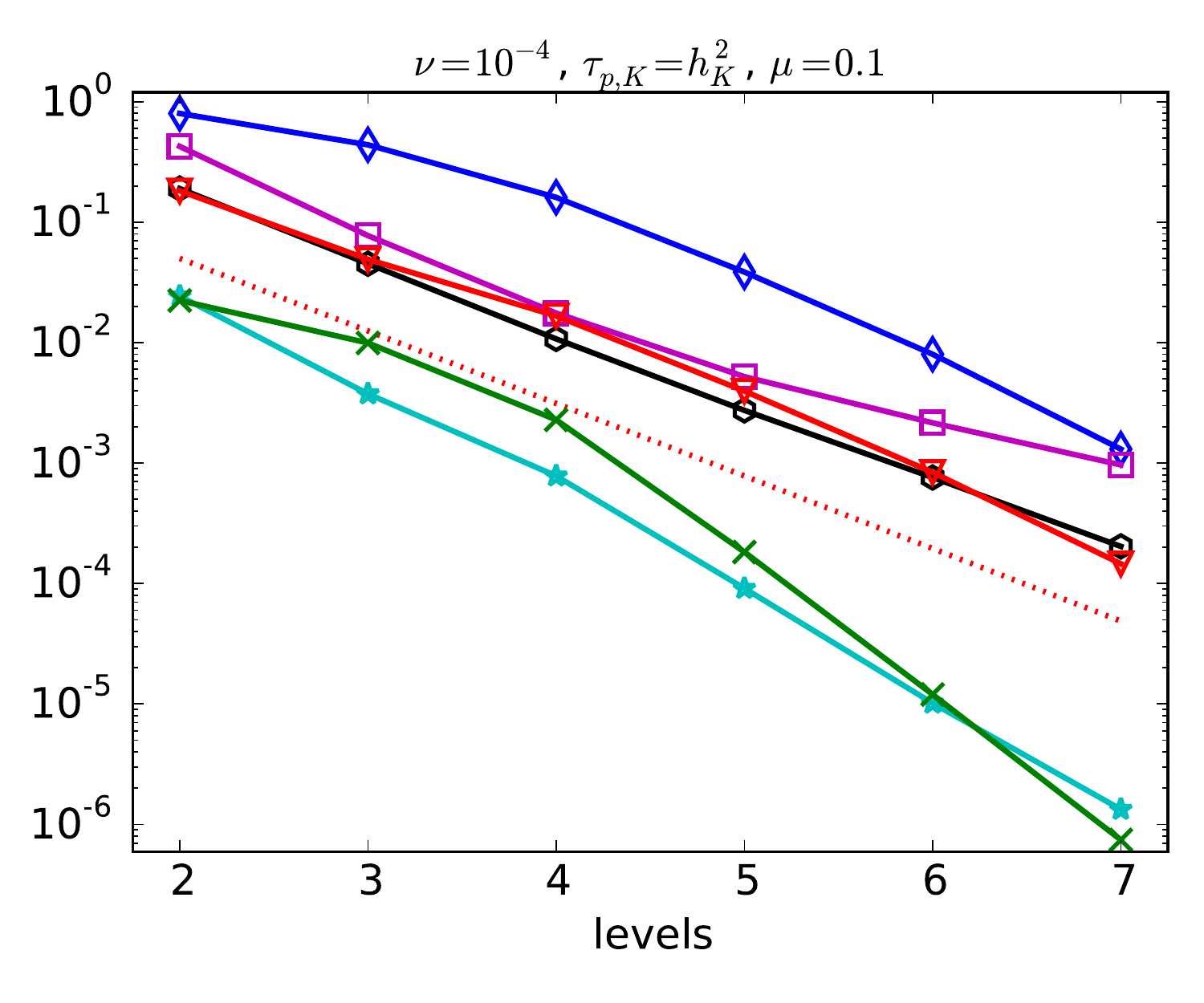}
\hspace*{1em}
\includegraphics[width=0.4\textwidth]{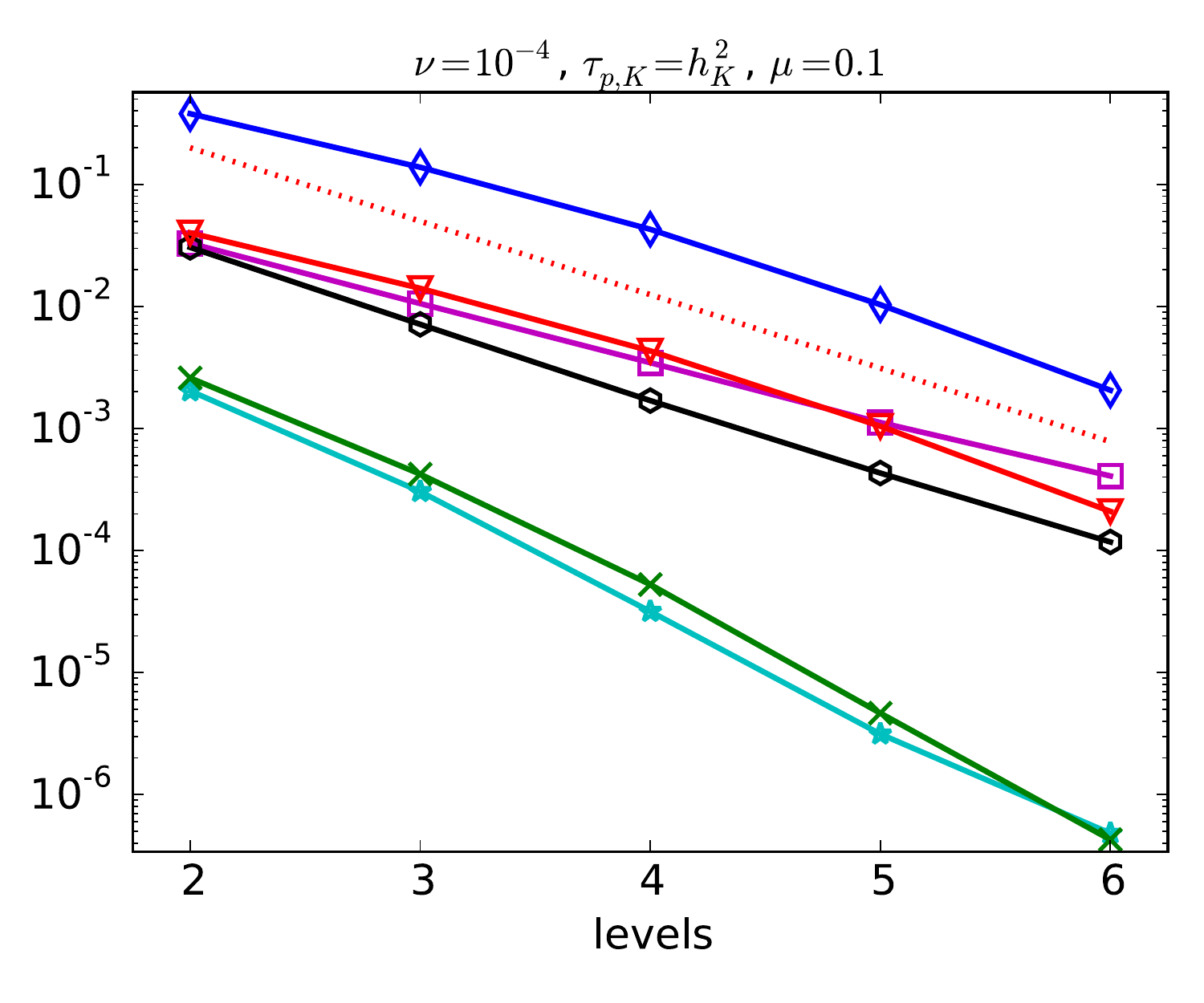}}
\centerline{\includegraphics[width=0.4\textwidth]
{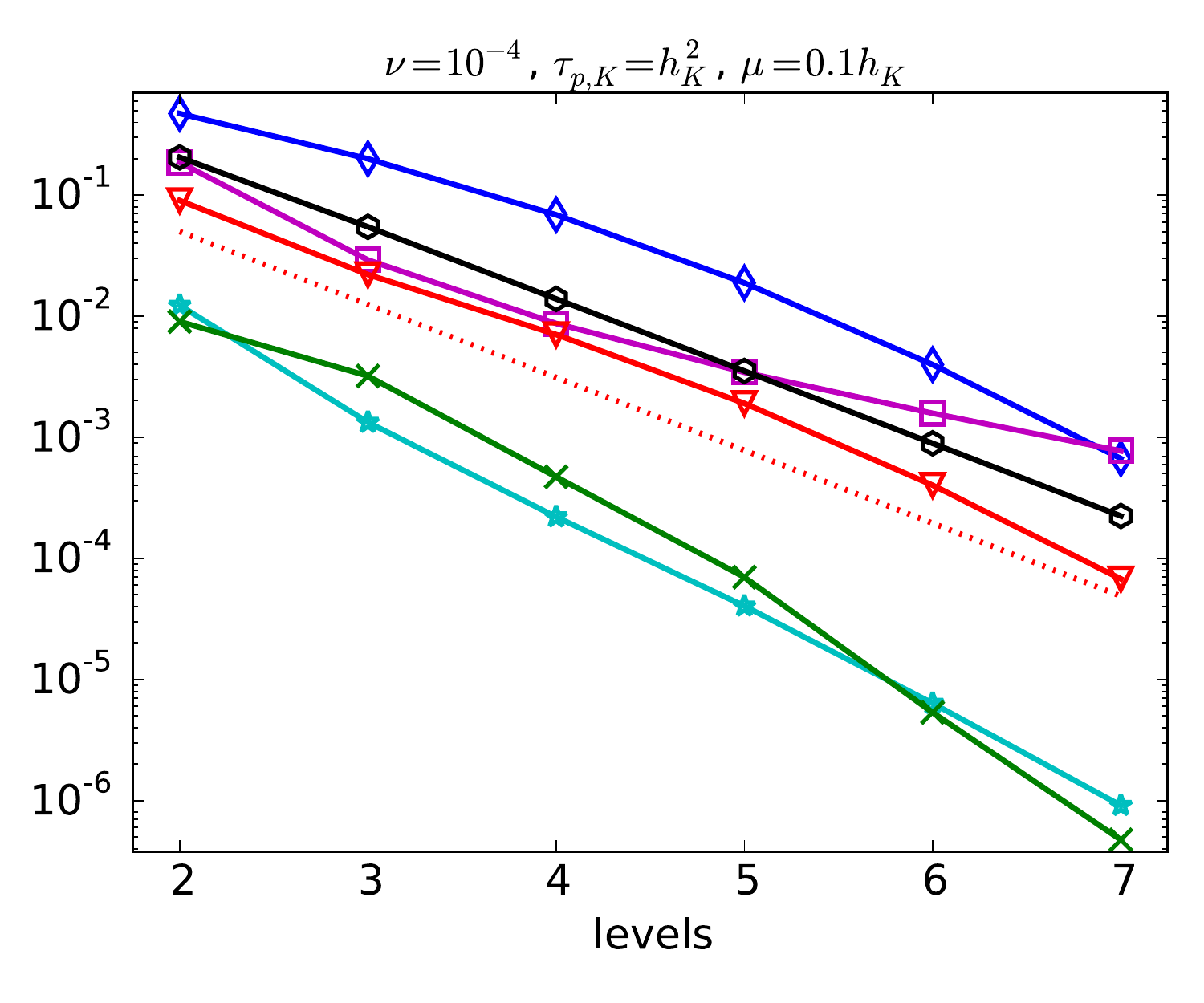}
\hspace*{1em}
\includegraphics[width=0.4\textwidth]{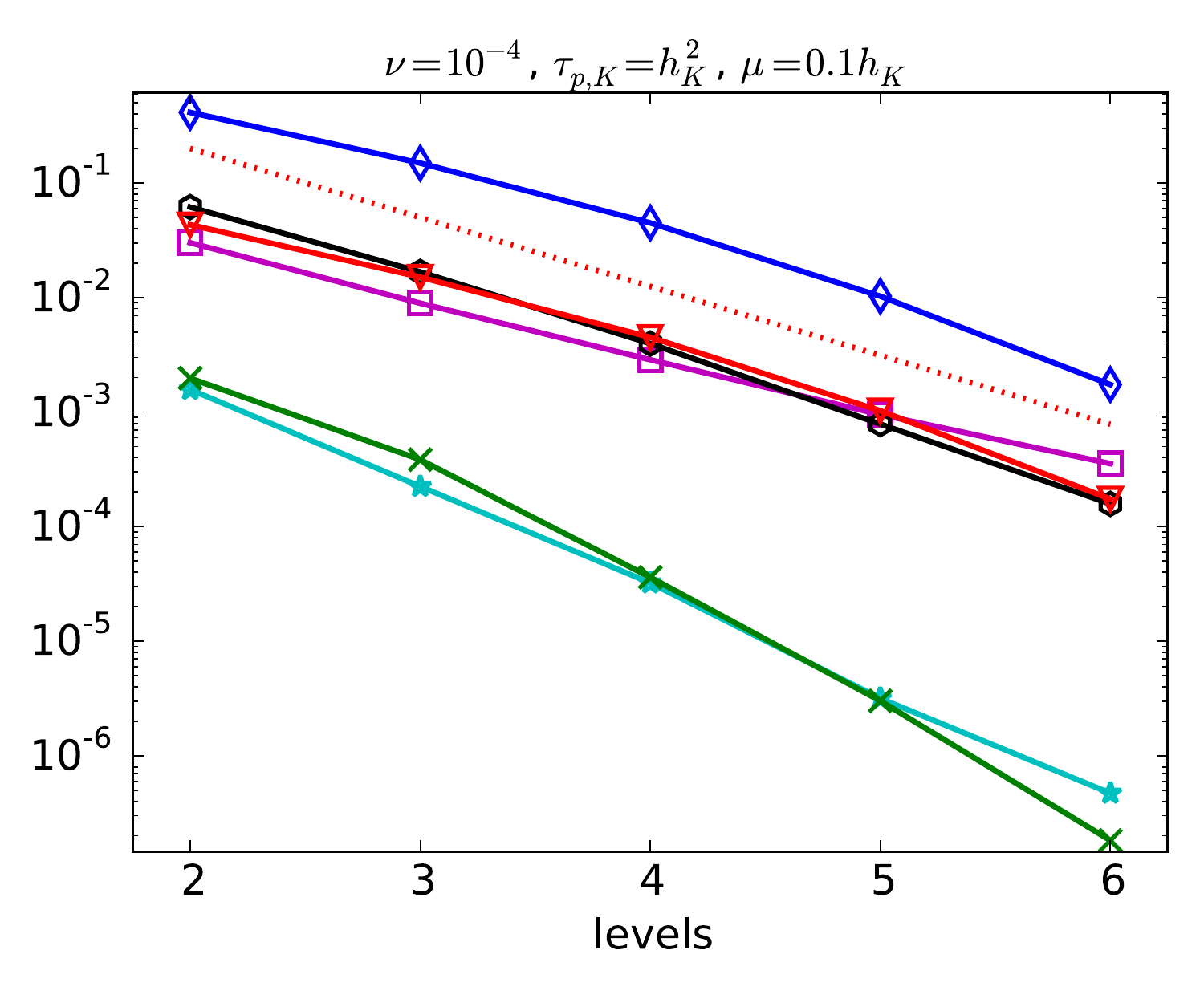}}
\centerline{\includegraphics[width=0.25\textwidth]
{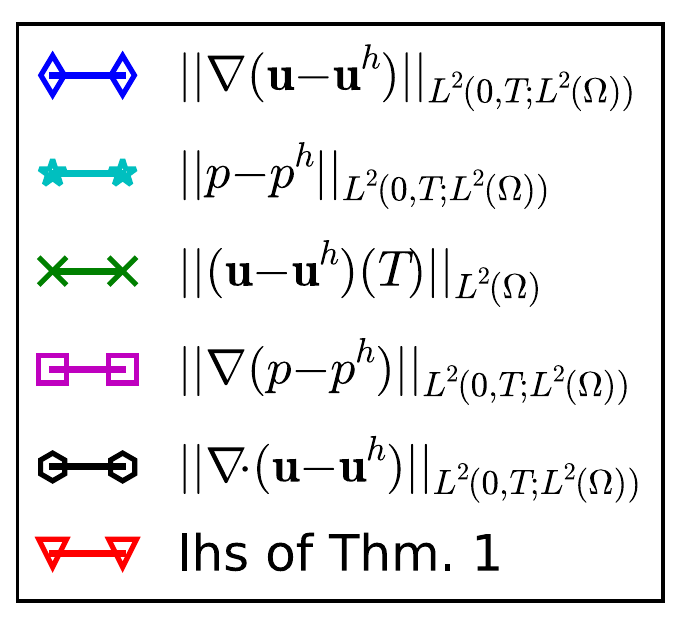}}
\caption{LPS with global grad-div stabilization, $P_2/P_2$ pair of finite element spaces, Grid~1 (left) and Grid~2 (right), dotted line: slope for second order convergence.}\label{fig:p2p2_glob_lps_glob_grad_div}
\end{figure}

\begin{figure}[t!]
\centerline{\includegraphics[width=0.4\textwidth]
{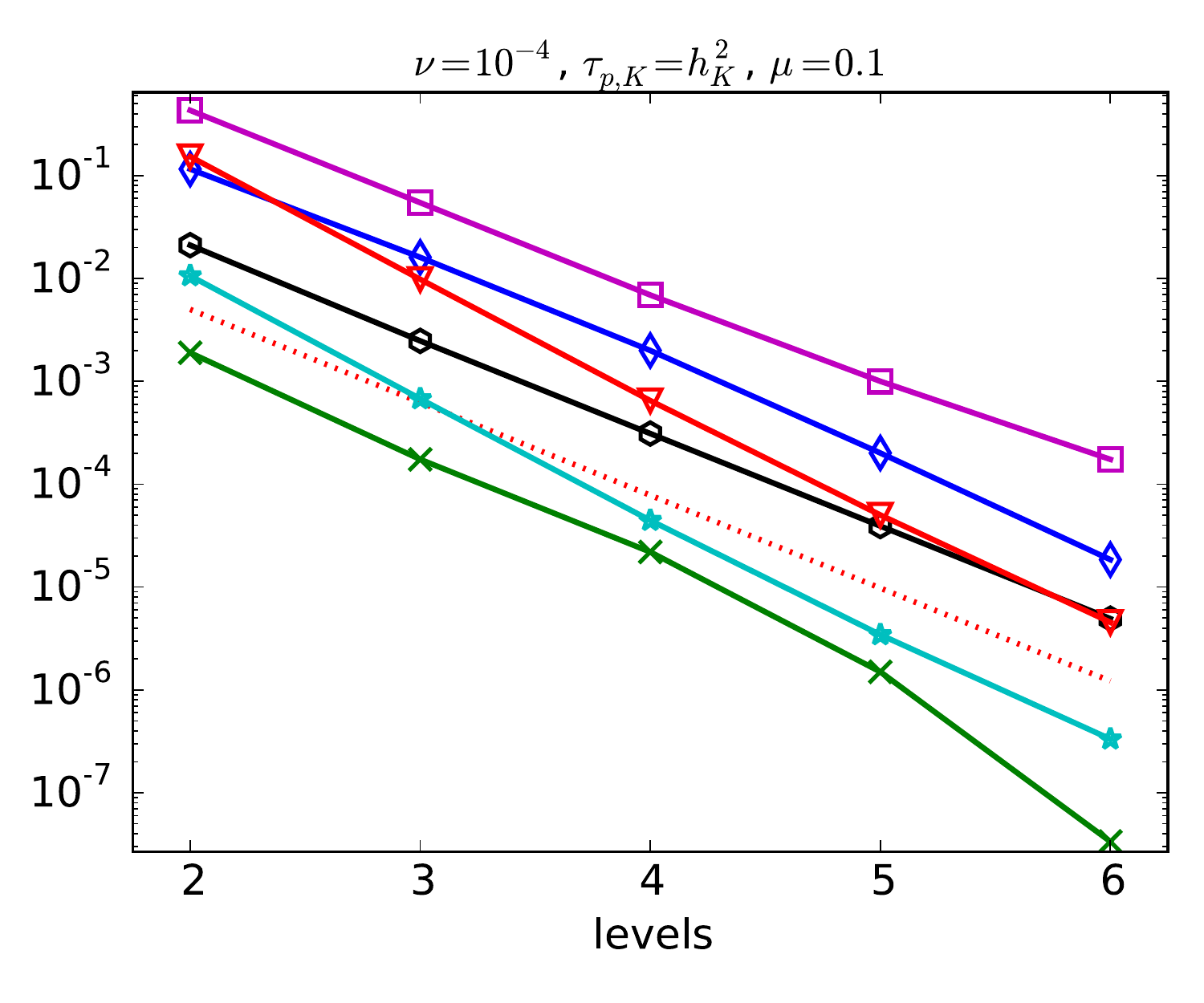}
\hspace*{1em}
\includegraphics[width=0.4\textwidth]{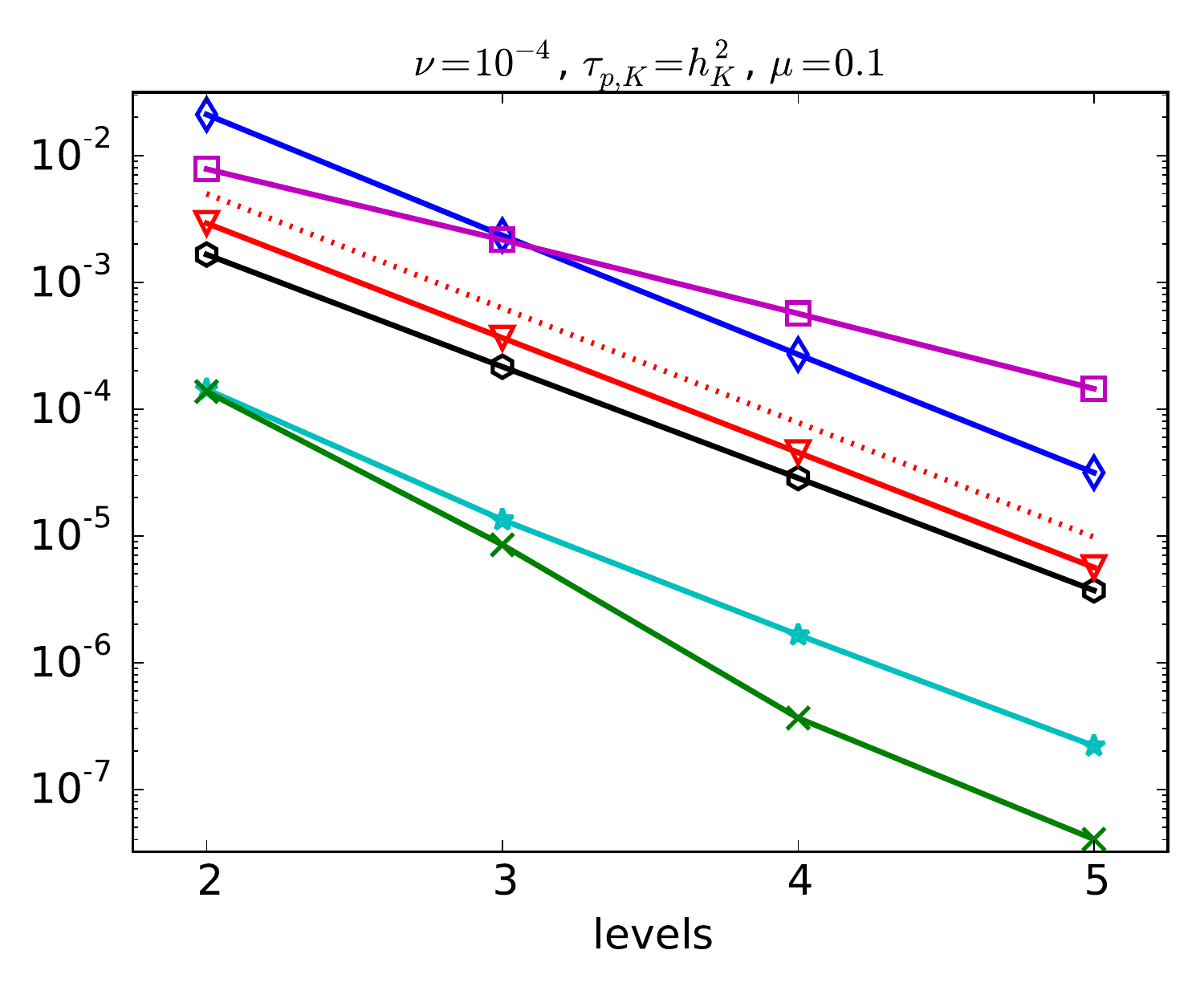}}
\centerline{\includegraphics[width=0.4\textwidth]
{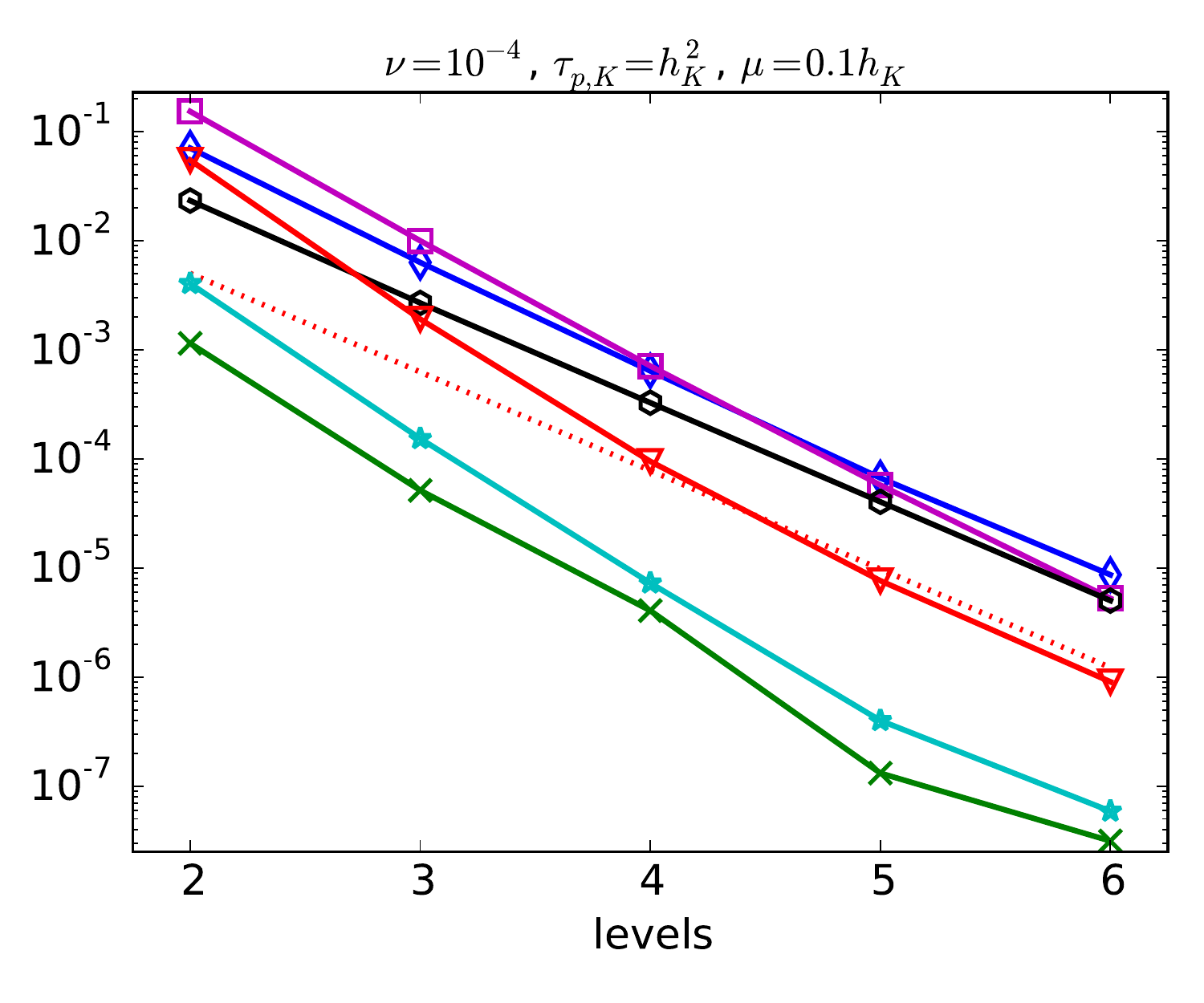}
\hspace*{1em}
\includegraphics[width=0.4\textwidth]{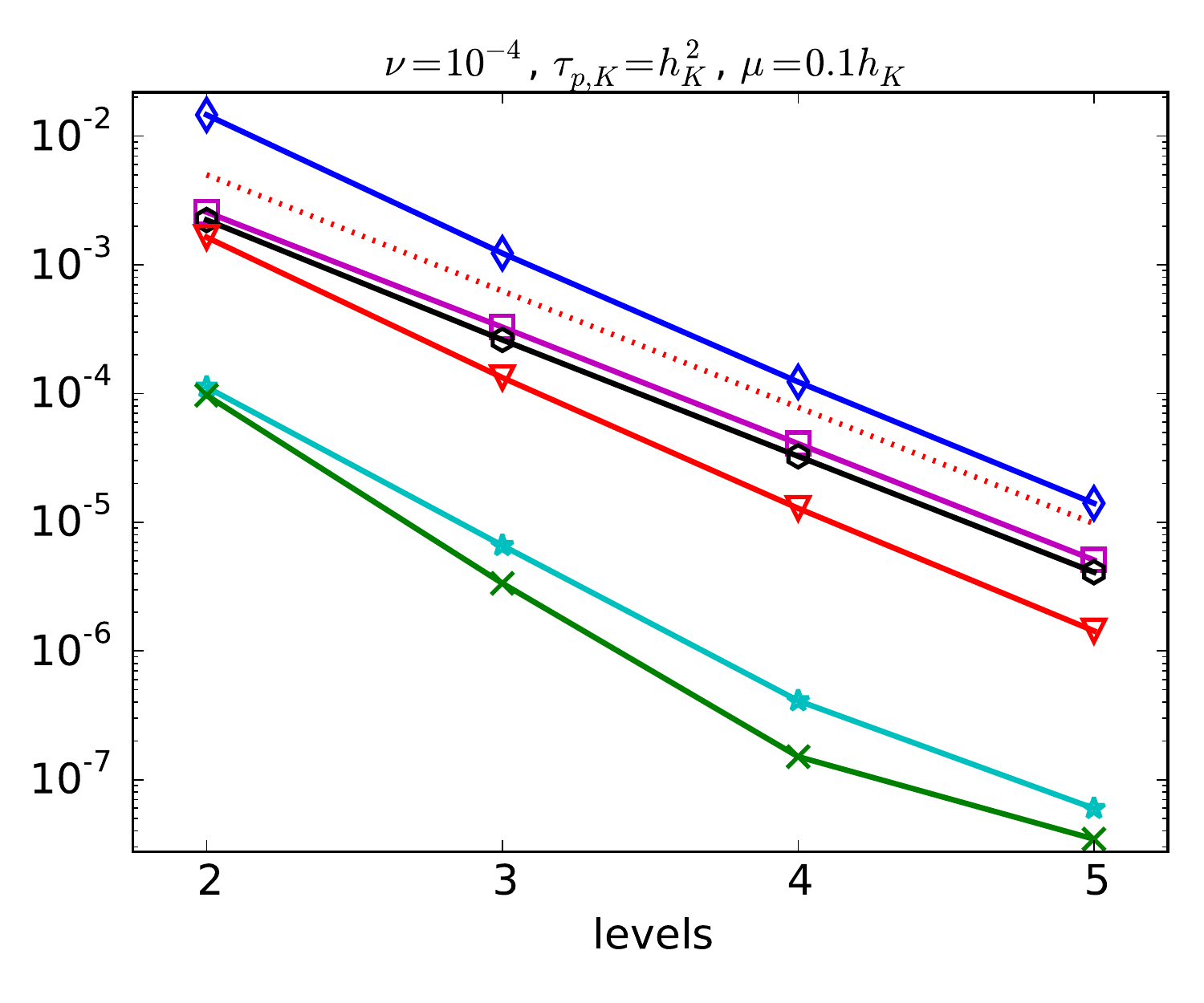}}
\caption{LPS with global grad-div stabilization, $P_3/P_3$ pair of finite element spaces,
Grid~1 (left) and Grid~2 (right), dotted line: slope for third order convergence, same legend as in Figure~\ref{fig:p2p2_glob_lps_glob_grad_div}.}\label{fig:p3p3_glob_lps_glob_grad_div}
\end{figure}

Results presented with the $P_2/P_2$ pair of finite elements are presented
in Figure~\ref{fig:p2p2_glob_lps_glob_grad_div} and with the $P_3/P_3$ pair
of spaces in Figure~\ref{fig:p3p3_glob_lps_glob_grad_div}. These results agree
with the analytical predictions. Concerning the grad-div stabilization parameter
there are only minor differences in the results. For the $P_3/P_3$ pair of spaces,
$\mu = 0.1h_K$ gives a somewhat better approximation of the pressure.

\begin{figure}[t!]
\centerline{\includegraphics[width=0.4\textwidth]
{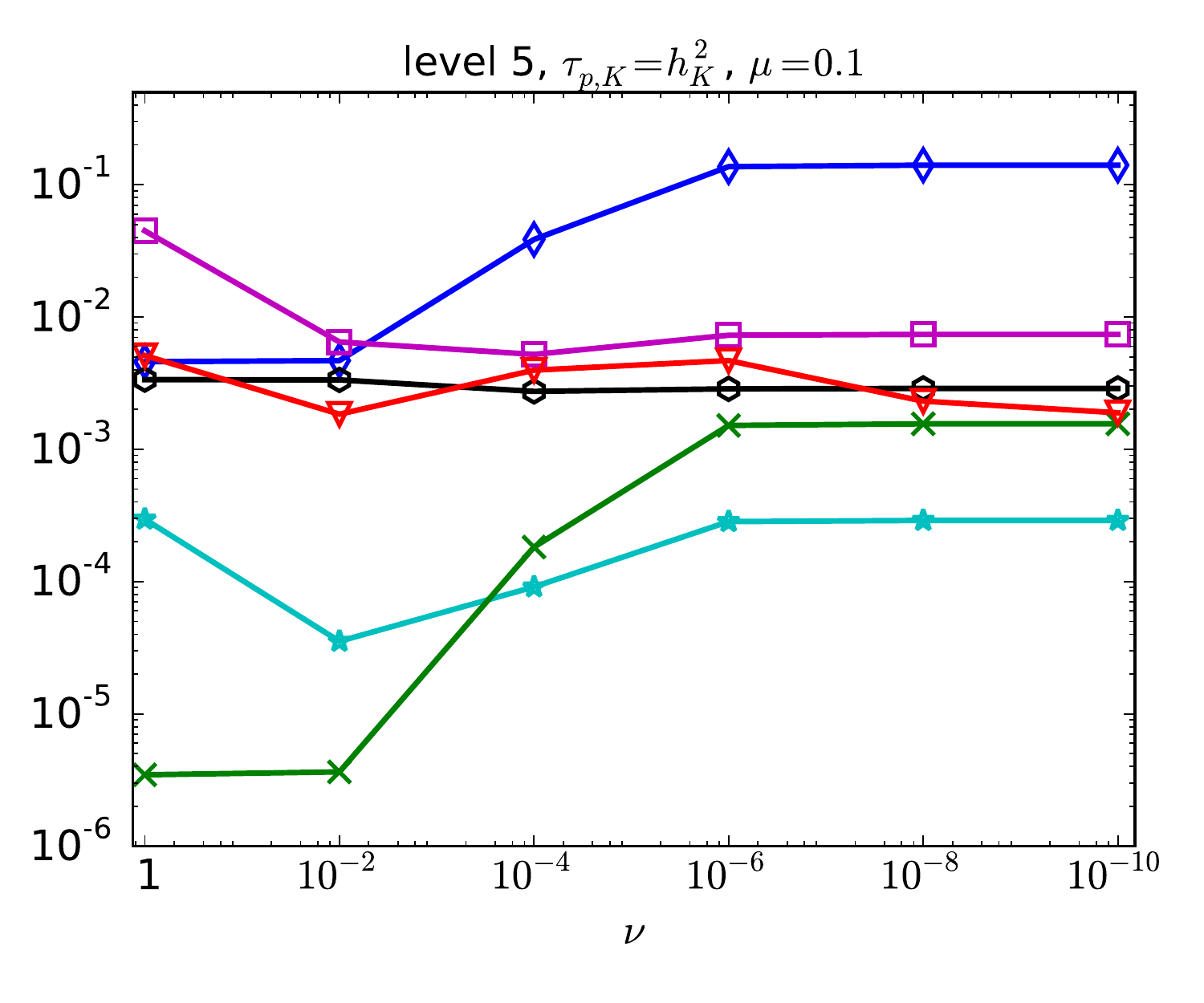}
\hspace*{1em}
\includegraphics[width=0.4\textwidth]{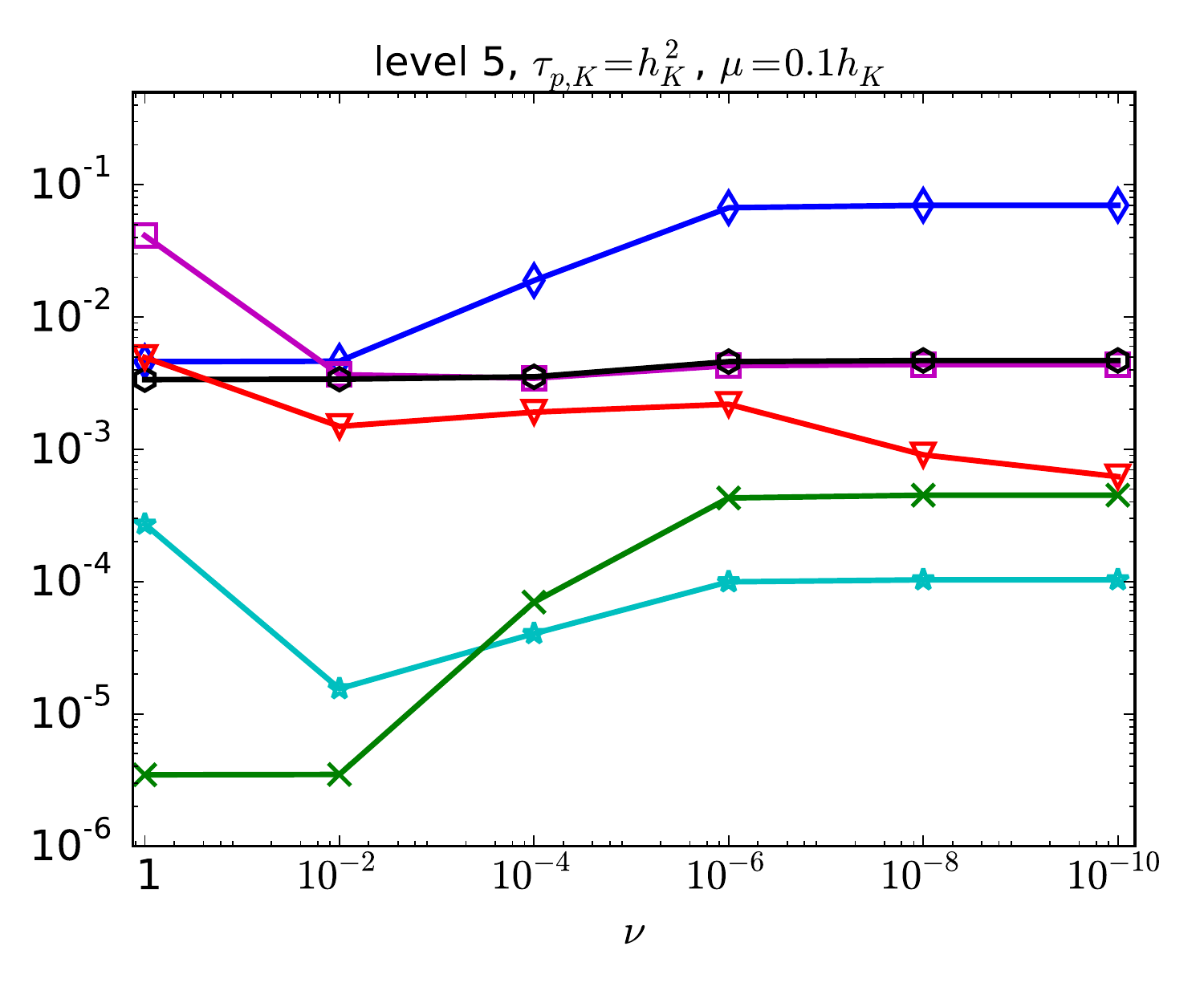}}
\caption{LPS with global grad-div stabilization, $P_2/P_2$ pair of finite element spaces,
Grid~1, behavior of errors with respect to $\nu^{-1}$, same legend as in Figure~\ref{fig:p2p2_glob_lps_glob_grad_div}.}\label{fig:fig:p2p2_glob_lps_glob_grad_div_lev5}
\end{figure}

Figure~\ref{fig:fig:p2p2_glob_lps_glob_grad_div_lev5} displays a representative result
for the dependency of the errors on the viscosity. It can be seen that all errors,
in particular the approximation of the error on the left-hand side of \eqref{eq:after_gron_2}, are bounded for $\nu \to 0$. This behavior coincides with the
analytical prediction.

\subsection{A method with rate of decay $s+1/2$  of the velocity error for $\nu \le h$}

Simulations for the method analyzed in Section~\ref{sec:ordersplusonehalf} were performed on
the irregular Grid~2, to prevent any superconvergence effects, for $\nu=10^{-8}$, such that
condition \eqref{eq:ass_sect6_1} is satisfied, and for the final time $T=0.5$. The remaining
setup of the simulations was as described in Section~\ref{sec:numres_method3}.

The methods incorporating the fluctuations of the velocity gradient were implemented as follows.
Generally, the nonlinear problems were solved with a fixed point iteration (Picard iteration).
Since the matrix representing the fluctuations of the gradient possesses a wider stencil
than all other matrices for the velocity-velocity coupling, we put the term with the
fluctuations of the velocity gradient on the right-hand side in the Picard iteration. In order
to achieve a satisfying rate of convergence of this iteration, numerical tests showed that the parameters
$\{\tau_{\nu,K}\}$ should be rather small. In addition, we could see that increasing these
parameters above a certain value leads to a notable increase of the errors. Altogether,
for the irregular Grid~2, $\tau_{\nu,K} = 0.01 h_K$ turned out to be an appropriate choice.
In view of condition \eqref{eq:cond_dif}, the LPS parameters for the pressure were chosen
to be $\tau_{p,K} = 10^{-4} h_K$.

\begin{figure}[t!]
\centerline{\includegraphics[width=0.4\textwidth]
{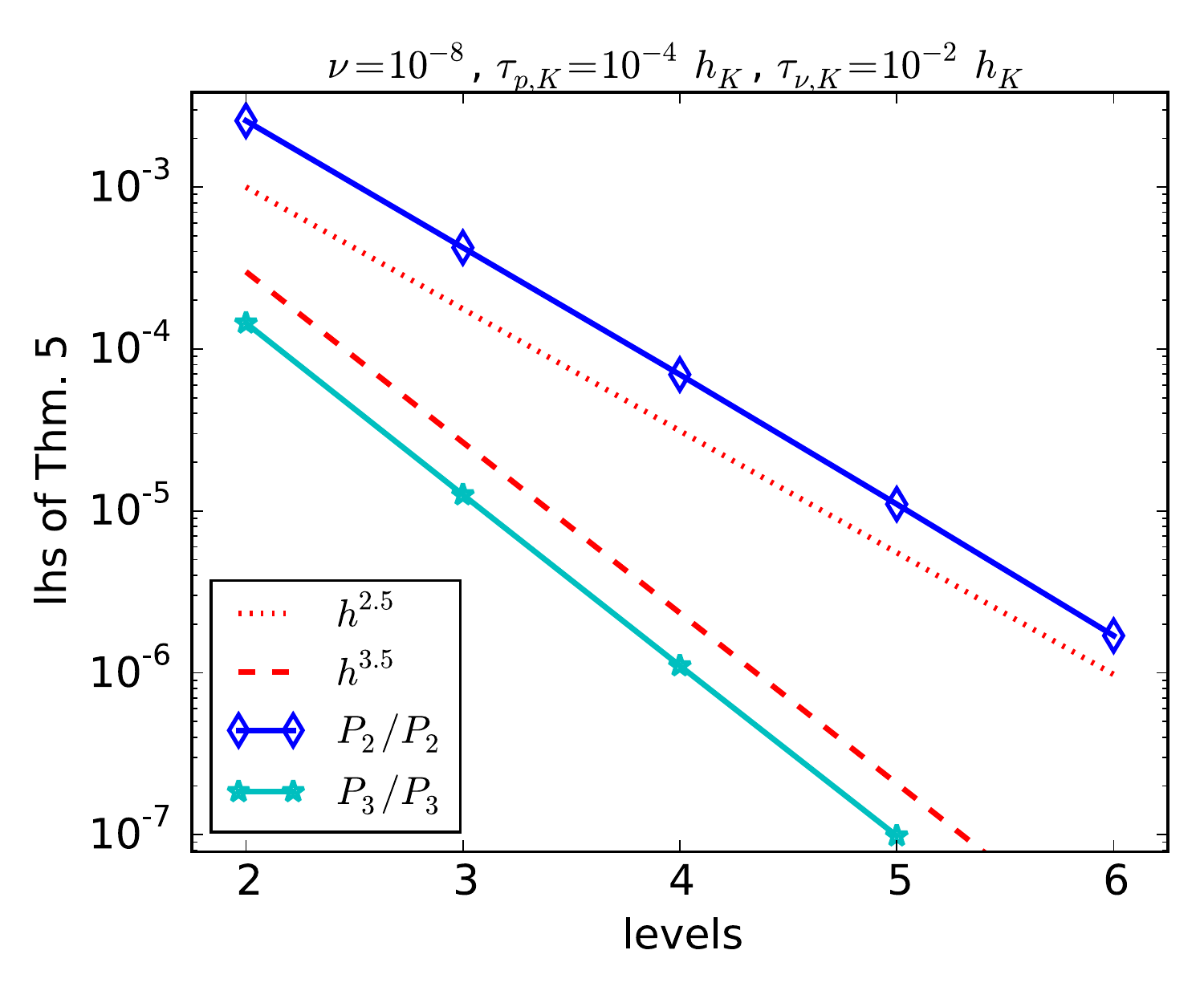}}
\caption{A method with rate of decay $s+1/2$  of the velocity error for $\nu \le h$,
computational results on Grid~2.}
\label{fig:numres_sect6}
\end{figure}

An error bound for the considered method was derived in Theorem~\ref{thm:velo_bound_4b}.
In the numerical simulations, the terms with the fluctuations on the left-hand side
of \eqref{eq:err_after_gron_22b} were
approximated by
$$
\Delta t\sum_{j=1}^n\|\sigma_h^*(\nabla (I_hp^j- p_h^j))\|_{\tau_p}^2, \quad
\Delta t\sum_{j=1}^n \|\sigma_h^*(\nabla (I_h\bu^j - \bu_h^{j}))\|_{\tau_{\nu}}^2,
$$
where $I_h$ is the Lagrangian interpolant. With the interpolants of the solution, these
terms can be simply computed by matrix-vector operations with the matrix of the fluctuations.

Computational results are presented in Figure~\ref{fig:numres_sect6}. One can observe the
proposed rates of decay of the velocity error. Having a detailed look on the individual contributions
of the error, we could see that the $L^2$ error and the
fluctuations of the velocity gradient were dominant.

\bibliographystyle{plain}

\end{document}